\documentclass[reqno]{amsart}
\numberwithin{equation}{section}

\usepackage{amsmath}
\usepackage{amsthm}
\usepackage{amssymb}
\usepackage{graphics}

\usepackage{mathrsfs} 

\usepackage{comment} 

\usepackage{graphicx}
\usepackage{amscd}
\usepackage{pstricks}
\usepackage{pst-plot}
\usepackage{multido}
\usepackage{pst-coil}
\usepackage{colortab}

\usepackage{comment} 

\usepackage{tikz} 

\usepackage{float} 

\usepackage[foot]{amsaddr} 

\makeatletter
\renewcommand{\email}[2][]{%
  \ifx\emails\@empty\relax\else{\g@addto@macro\emails{,\space}}\fi%
  \@ifnotempty{#1}{\g@addto@macro\emails{\textrm{(#1)}\space}}%
  \g@addto@macro\emails{#2}%
}
\makeatother

\theoremstyle{plain}
\newtheorem{theorem}{\bf Theorem}[section]
\newtheorem{proposition}[theorem]{\bf Proposition}
\newtheorem{lemma}[theorem]{\bf Lemma}
\newtheorem{corollary}[theorem]{\bf Corollary}

\theoremstyle{remark}
\newtheorem{remark}[theorem]{\bf Remark}

\begin{document}
\title{A priori estimates of Mizohata-Takeuchi type for the Navier-Lam\'e operator}
\author{J.A. Barcel\'o$^1$}
\author{A. Ruiz$^2$}
\author{M.C. Vilela$^1$}
\author{J. Wright$^4$}
\address{$^1$Universidad Politécnica de Madrid, Spain.}
\address{$^2$Universidad Autónoma de Madrid, Spain.}
\address{$^4$University of Edingburgh, UK}
\email{juanantonio.barcelo@upm.es, alberto.ruiz@uam.es}
\email{maricruz.vilela@upm.es and J.R.Wright@ed.ac.uk}
\thanks{The first and the third authors were supported by the Spanish Grant PID2021-124195NB-C31, the second one by the the Spanish Grant PID2021-124195NB-C31 and the fourth by a Leverhulme Research Fellowship RF-2023-709$\backslash$9.}
\subjclass[2010]{Primary 35J47, 74B05. Secondary 42B37.}
\keywords{Elasticity system, a priori estimates.}
\begin{abstract}
The Mizohata-Takeuchi conjecture for the resolvent of the Navier-Lamé equation is a weighted estimate with weights in the so-called Mizohata-Takeuchi class for this operator when one approaches the spectrum (Limiting Absorption Principles).
We prove this conjecture in dimensions 2 and 3 for weights with a radial majorant in the Mizohata-Takeuchi class.
This result can be seen as an extension of the analogue for the Laplacian given in \cite{BRV1997}.
We also prove that radial weights in this class are not invariant for the Hardy-Littlewood  maximal function, hence the methods in \cite{BFPRV2012} used to extend estimates for the Laplacian to the Navier-Lamé case, do not work.
\end{abstract}
\maketitle
\markboth{J.A. Barcel\'o, R. Ruiz, M.C. Vilela and J. Wright}{A priori estimates of Mizohata-Takeuchi type for the Navier-Lam\'e operator}
\section{Introduction and statement of results}
Nowadays, the Mizohata-Takeuchi conjecture is known as a weighted estimate for the extension operator of the Fourier transform on the sphere.
More precisely, if $S_R^{d-1}$ denotes the sphere of radius $R$ in $\mathbb{R}^d$, $d\sigma_R$ the surface measure induced by the Lebesgue measure, $V:\mathbb{R}^d\longrightarrow [0,+\infty)$ and
$f\in L^2(S_R^{d-1})$, the conjecture can be formulated as the following weighted estimate for the extension operator over the sphere:
			\begin{equation}
   \label{conjetura}
				\left\|\widehat{fd\sigma_R}\right\|^2_{L^2(V)}
				\le
				c\,
				\textrm{sup}\left\{\int_{\textrm{L}}V:\textrm{L} \text{ line in }\mathbb{R}^d  \right\}\,
				\|f\|^2_{L^2(S_R^{d-1})}.
			\end{equation}
The integral of $V$ along the line $\textrm{L}$ is known as the X-ray transform of $V$, this is the function $XV$ defined on all the lines in $\mathbb{R}^d$ by
   \begin{equation}
       \label{transformadaRX}
       XV(\textrm{L}):=\int_{\textrm{L}}V=\int_{-\infty}^\infty V(y+t\theta)dt,
   \end{equation}
where $y\in\mathbb{R}^d$ and $\theta\in S^{d-1}$ determine the line $\textrm{L}$.

In the  case of radial weights $V$,  \eqref {conjetura} was proved in \cite{BRV1997} and \cite{CS1997} independently.
The Mizohata-Takeuchi conjecture has been an object of interest related to Stein's conjecture in Fourier  Analysis, see \cite{S1978}, \cite{BBC2008}, \cite{Sh2022}, \cite{BN2021}, \cite{BS2019}, \cite{CIW2023} and the references therein.

The origin of this conjecture is an initial value problem. 
In the 80's, Takeuchi (see \cite{T1974},\cite{T1980} and \cite{T1985}) was studying the following initial value problem associated to the linear Schrödinger equation perturbed by its first order term:
\begin{equation}
		\label{ivp}
		\left\{
		\begin{array}{ll}
		i^{-1}\partial_t u+\Delta_xu+\mathbf{b}(x)\cdot\nabla_x u=F(x,t),&\quad x\in\mathbb{R}^d, t\in\mathbb{R},
		\\
		u(x,0)=u_0(x).
		\end{array}
		\right.\
		\end{equation}
with $u_0\in L^2(\mathbb{R}^d)$, $F\in C(\mathbb{R};L^2(\mathbb{R}^d))$ complex-valued functions, and 
$\mathbf{b}$ a vector-valued function.
He claimed that \eqref{ivp} is well-posed in $L^2(\mathbb{R}^d)$ if and only if
    \begin{equation}
    \label{Takeuchi}
	\textrm{sup}\left\{\textrm{Re}\int_0^\infty\mathbf{b}(x+s\theta)\cdot\theta\, ds: x\in\mathbb{R}^d, \theta\in S^{d-1} \right\}<\infty.
    \end{equation}
Later Mizohata (see \cite{M1985}) noticed that Takeuchi's proof was not clear and showed that condition \eqref{Takeuchi} is necessary but left open the question of whether it is also a sufficient condition. 
As far as we know, this Mizohata-Takeuchi conjecture is not yet resolved.

In 1997, Barceló-Ruiz-Vega (see \cite{BRV1997}) introduced a scalar version of \eqref{Takeuchi}. 
More precisely, they consider scalar potentials $V:\mathbb{R}^d\longrightarrow [0,+\infty)$ satisfying the \emph{scalar Mizohata-Takeuchi} condition
	\begin{equation}
	    \label{escalarTakeuchi}
     |||V|||_X:=\textrm{sup}\left\{\int_{-\infty}^{\infty}V(x+s\theta)\, ds: x\in\mathbb{R}^d,\theta\in S^{d-1} \right\}<\infty.
	\end{equation}
We will call the supremum appearing in \eqref{escalarTakeuchi} the X-ray norm of $V$ and we will denote the class of potentials satisfying \eqref{escalarTakeuchi} by $\mathcal{T}(\mathbb{R}^d)$.
They also introduced the \emph{radial Mizohata-Takeuchi} class, denoted by $\mathcal{T}_{rad}(\mathbb{R}^d)$, as the class of
nonnegative functions $V$ such that the function
\begin{equation}
\label{mayorante_radial}
    \widetilde V(r):=\sup_{\omega\in S_r^{d-1}} V(\omega)
\end{equation}
satisfies
\begin{equation}
    \label{normaTakeuchi}
    |||V|||:=\sup_{\mu\ge0}\int_{\mu}^{\infty}\frac{r\,\widetilde V(r)}{(r^2-\mu^2)^{1/2}}\,dr<\infty.
\end{equation}

Abusing notation, we write $\widetilde V(x)=\widetilde V(r)$ where $r=|x|$. One can see that condition \eqref{normaTakeuchi} means that the $X-$ray transform of the radial majorant $\widetilde V(x)$ is bounded everywhere. 
In fact, we may assume that $x\cdot \theta=0$ in \eqref{escalarTakeuchi} and then, by making the change of variable $x+s\theta=r$ in the integral with $V=\widetilde V$ we obtain the integral in \eqref{normaTakeuchi}.

In \cite[Theorem 4]{BRV1997} they proved that the initial value problem \eqref{ivp} is well-posed in the homogeneous Sobolev space $\overset{\cdot}{H}^{1/2}(\mathbb{R}^d)$
when $\mathbf{b}$ is bounded by a scalar radial function $V$ with a sufficiently small X-ray norm.

The essential tool in their proof is an a priori weighted estimate for the solution of the stationary Helmholtz equation
\begin{equation}
    \label{Helmholtz_eq}
    \Delta u(x)+k^2 u(x)=f(x),\qquad k>0,\,x\in\mathbb{R}^d,
\end{equation}
satisfying the Sommerfeld outgoing radiation condition
\begin{equation}
    \label{SRC}
    (\partial_r-ik)u=o(r^{-(d-1)/2}),\qquad
r=|x|\rightarrow\infty.
\end{equation}
They proved that if $f\in C_0^\infty(\mathbb{R}^d)$ and $V$ is a radial weight with a bounded X-ray norm, then the solution $u$ of \eqref{Helmholtz_eq}-\eqref{SRC} satisfies the following a priori estimate (see \cite[Th. 1 and Th. 2]{BRV1997}):
    \begin{equation}
    \label{l-resolvent}
	\|\nabla u\|_{L^2(V)}+k\|u\|_{L^2(V)}
	\le
	c\,
	|||V|||\,
	\left\|f\right\|_{L^2(V^{-1})}.
    \end{equation}
By a limiting absorption principle argument, they got the same estimate for
the outgoing resolvent of the Laplace operator. This resolvent can be written as the following sum:
\begin{equation*}
		\label{resolvente_Lapace}
		(\Delta+k^2+i0)^{-1}f(x)
		=\frac{i}{k}\,\widehat{d\sigma_k}\ast f(x)
		+\textrm{p.v.}\int_{\mathbb{R}^d}\frac{\widehat{f}(\xi)}{|\xi|^2-k^2}\,e^{ix\cdot\xi}d\xi.
\end{equation*}
So, formally, we can say that taking the imaginary part in the estimate for the resolvent,
they obtained the same estimate for the convolution Stein-Tomas operator.
In fact, they proved a stronger result (see \cite[Theorem 3]{BRV1997}). More precisely, they characterized the radial weights for which the following estimate holds:
\begin{equation*}
    \left\|\widehat{d\sigma_k}\ast f\right\|_{L^2(V)}
    \le
    c\,
    |||V|||\,
    \|f\|_{L^2(V^{-1})}.	
\end{equation*}
And also the equivalent estimates for the extension and restriction operator on the sphere, which are related to \eqref{conjetura}.
    
The purpose of this article is to study this problem in the context of linear elasticity, which 
is governed by the equation
\begin{equation}
\Delta^{\ast}\mathbf{u}(x)+\omega^{2}\mathbf{u}(x) = \mathbf{f}(x), \qquad\qquad \omega>0,\ x\in\mathbb{R}^{d}, 
\label{ecuacion}
\end{equation}
where $\mathbf{u},$ the displacement vector, is a vector-valued function from
$\mathbb{R}^{d}$ to $\mathbb{R}^{d}$ and,
\begin{equation}
\label{operador}\Delta^{\ast}\mathbf{u}(x) = \mu\Delta\mathbf{u}(x)+(\lambda+\mu)\nabla div\,\mathbf{u}(x),
\end{equation}
with $\Delta\mathbf{u}$ defined componentwise. 
Here $\omega>0$ denotes the frequency of the time-harmonic wave, the constants $\lambda$ and $\mu$ are known as the Lam\'{e} constants and $\Delta^*$ is  the \emph{Navier-Lam\'e operator}. 

Throughout this paper we will assume that $\mu>0$ and
$2\mu+\lambda>0$ so that the operator $\Delta^{\ast}$ is strongly
elliptic and, we will denote by $k_p$ and $k_s$ the speed of propagation of  longitudinal and  transversal waves, respectively, where
\begin{equation}
\label{kpks} \displaystyle k_p^2=\frac{\omega^2}{(2\mu+\lambda)}
\qquad \texttt{\rm{and}} \qquad k_s^2=\frac{\omega^2}{\mu}.
\end{equation}

It is well-known that any solution $\mathbf{u}$ of the \emph{homogeneous Navier-Lam\'e equation}
\begin{equation}
\label{homogenea}
\Delta^{\ast}\mathbf{u}(x)+\omega^{2}\mathbf{u}(x)=\mathbf{0},
\end{equation}
with $\Delta^{\ast}$ given by \eqref{operador}, in a domain, can
be written as the sum of the so called compressional part, denoted
by $\mathbf{u}_p,$ and the shear part, denoted by $\mathbf{u}_s,$
where
\begin{equation}
\label{upus} \displaystyle \mathbf{u}_p=-\frac{1}{k_p^2}\nabla
div\, \mathbf{u} \qquad \texttt{\rm{and}} \qquad
\mathbf{u}_s=\mathbf{u}-\mathbf{u}_p.
\end{equation}
Observe that $\mathbf{u}_p$ and $\mathbf{u}_s$ are solutions of
the \emph{vectorial homogeneous Helmholtz equations} $\Delta
\mathbf{u_p}(x)+k_p^{2}\mathbf{u_p}(x)=\mathbf{0}$ and $\Delta
\mathbf{u_s}(x)+k_s^{2}\mathbf{u_s}(x)=\mathbf{0},$ respectively.

Besides, if $\mathbf{u}$ is an entire solution (i.e. a solution in
the whole $\mathbb{R}^d$) of \eqref{homogenea} satisfying the outgoing
Kupradze radiation conditions:
\begin{eqnarray}
\label{radiacionup}
(\partial_r-ik_p)\mathbf{u}_p&=&\mathbf{o}(r^{-(d-1)/2}),\qquad
r=|x|\rightarrow\infty,
\\
\label{radiacionus}
(\partial_r-ik_s)\mathbf{u}_s&=&\mathbf{o}(r^{-(d-1)/2}),\qquad
r=|x|\rightarrow\infty,
\end{eqnarray}
then, $\mathbf{u}=\mathbf{0}$ (see \cite{Kupradze} for the
three-dimensional case). 

As a consequence, for a vector-valued function 
$\mathbf{f}\in\mathcal{C}_0^{\infty}(\mathbb{R}^d),$ if there exists a solution of the
\emph{Navier-Lam\'e equation} \eqref{ecuacion}
satisfying the Kupradze radiation conditions \eqref{radiacionup}
and \eqref{radiacionus}, where $\mathbf{u}_p$ and $\mathbf{u}_s$
are given by \eqref{upus} outside the support of $\mathbf{f},$ then
the solution is unique.


In \cite{BFPRV2012}    several classes  of weights $V$ are considered, those for  which  estimate \eqref{l-resolvent} holds, and  the analogous  estimate was extended    to the resolvent  of   Navier-Lam\'e  operator in elasticity.   The  extension was based on the good behavior of the different classes  of  weights with respect to  singular integral operators, Leray's decomposition of the  resolvent of Navier operator in terms  of the resolvent  of Laplacians at different energies (longitudinal and transversal  waves), together with Rubio de Francia  algorithm, (see \cite{GR1985}).
Nevertheless the case of the radial Mizohata-Takeuchi class  remained an open problem. The question was the invariance of this class under Riesz  transforms (see Remark 3.11 in \cite{BFPRV2012}).

In this work we give a negative answer to this question by giving a bounded sequence of radial functions in $\mathcal{T}(\mathbb{R}^d)$ whose image by the Hardy-Littlewood  maximal operator is unbounded (see Proposition \ref{Contraejemplo}).
Nevertheless, we are able to extend estimate \eqref{l-resolvent} to the Navier-Lamé operator for radial potentials (see Theorem \ref{principal}). 
This estimate provides weighted estimates for the resolvent (see Corollary \ref{cor_principal}) which are an important tool to solve boundary and  scattering inverse problems  for perturbation of the  free operator, see \cite{BFPRV2018} for Navier-Lam\'e  operator.

The main results of this paper are the following.
\begin{theorem}
	\label{principal}
Let $V\in \mathcal{T}_{rad}(\mathbb{R}^d)$ with $d=2,3$.
If $\mathbf{f}\in \mathcal{C}_0^{\infty}(\mathbb{R}^d),$ then there exists $c>0$ independent of $\mathbf{f}$, $\omega$ and $V$ such that the following a priori estimate for $\mathbf{u}$, the unique solution of \eqref{ecuacion} satisfying \eqref{radiacionup} and \eqref{radiacionus}, holds:
	\begin{equation}
	\label{estimacion_u}
	\int_{\mathbb{R}^d}|\mathbf{u}(x)|^2V(x)dx
	\le \frac{c}{\omega^2}|||V|||^2
	\int_{\mathbb{R}^d}|\mathbf{f}(x)|^2V^{-1}(x)dx.
	\end{equation}
\end{theorem}
\begin{theorem}
	\label{teorema_gradiente}
Under the conditions of Theorem \ref{principal},
writing $\mathbf{u}=(u_1,\ldots,u_d)$, the following a priori estimate holds:
 \begin{equation}
	\label{estimacion_gradu}
 \sup_{1\le j\le d}\int_{\mathbb{R}^d}|\nabla u_j(x)|^2V(x)dx
	\leq c\,|||V|||^2
	\int_{\mathbb{R}^d}|\mathbf{f}(x)|^2V^{-1}(x)dx.
	\end{equation}
\end{theorem}
We would like to point out that Theorems \ref{principal} and \ref{teorema_gradiente} can be stated with two different weights $V_1$ and $V_2$ (see Theorem \ref{solucion_th} for the Helmholtz case).

As we mentioned before, in the study of well-posedness of initial value problems, a priori estimates for the associated stationary  problem are essential tools (see \cite{BRV1997} and \cite{BBR2006} for the Helmholtz case with Mizohata-Takeuchi weights). So estimates \eqref{estimacion_u} and \eqref{estimacion_gradu} can be used to prove existence and uniqueness of solutions of initial value problems associated to the evolution Navier-Lamé equation. And one obtains a priori estimates for these solutions (see \cite[Theorems 1.6 and 1.7]{BFPRV2012}). We will treat these issues in a forthcoming paper.

As a consequence of Theorem \ref{principal}, by using standard techniques (see \cite{BRV1997}) we get the following corollary.
\begin{corollary}
	\label{cor_principal}
Let $V\in \mathcal{T}_{rad}(\mathbb{R}^d)$ with $d=2,3$.
\begin{itemize}
	\item [(i)]
	If $z=\gamma+i\varepsilon$ with $\varepsilon\neq 0$, then there exists $c>0$ independent of $z$ and $V$ such that for any $\mathbf{u}\in\mathcal{C}_0^{\infty}(\mathbb{R}^d),$ the following a priori estimate holds:
	\begin{equation}
	\label{estimacion_resolvente}
	\int_{\mathbb{R}^d}|\mathbf{u}(x)|^2V(x)dx
	\le |z|^{-1}|||V|||^2
	\int_{\mathbb{R}^d}|(\Delta^{\ast}+z\mathrm{I})\mathbf{u}(x)|^2V^{-1}(x)dx.
	\end{equation}
	\item [(ii)] (Weak limiting absorption principle.)
    Let us denote by $\mathbf{R}(z)$ the extension to $\mathbf{L}^2(V^{-1}dx)$ of the operator 
    $(\Delta^{\ast}+z\mathrm{I})^{-1}$ defined in 
	$\mathbf{L}^2(\mathbb{R}^d)\cap\mathbf{L}^2(V^{-1}dx)$. Then for $\mathbf{f}\in \mathbf{L}^2(V^{-1}dx)$, the weak-limit
    \begin{equation*}
    \label{nweak}
    \displaystyle \mathbf {R}(\omega^2+i0)\mathbf {f}:=
    \underset{z\rightarrow\omega^2,\ \Im z>0}{\rm weak-lim}\ \mathbf {R}(z)\mathbf {f}
    \end{equation*}
    exists in $\mathbf{L}^2(Vdx)$ and is a weak solution of
    \eqref{ecuacion}. If $\mathbf{f}\in\mathcal{C}_0^{\infty}(\mathbb{R}^d),$ then 
    $\mathbf {R}(\omega^2+i0)\mathbf {f}$ is the solution satisfying the outgoing Kupradze radiation conditions \eqref{radiacionup} and \eqref{radiacionus}.
\end{itemize}
\end{corollary}

The structure of the paper is as follows. 
In  Section 2 we  prove that the radial Mizohata-Takeuchi class is not invariant under the  centered  Hardy-Littlewood  maximal function.  This  was the main obstacle  to extend the estimate  for  Helmoltz  equation to the Navier-Lam\'e equation in \cite{BFPRV2012}.

In Section 3 we give an addition formula  for the fundamental  solution  of the eigenfunction equation of Lam\'e-Navier operator. We  give  the proof in dimension  $d=2$ and refer  to \cite{Addition_Lame_3D} for  $d=3$. 
These formulae are a fundamental part in the proof of estimate \eqref{estimacion_u}.
We think that these formulae could be useful  in  other contexts, see  \cite{BCMM 2021}.

In Section 4 we prove Theorem \ref{principal}, up to some lemmas which we establish in  Section 6.
For simplicity, we give a detailed proof for $d=2$ and a sketch for $d=3$.
The proof of Theorem \ref{teorema_gradiente} is given in Section 5.

The proof  of  Theorem \ref{principal}, although it  follows the lines of the case of the Helmholtz equation in \cite{BRV1997}, it is much more involved. We will need new, delicate cancellations between longitudinal and transversal projections of solutions of the Navier-Lamé equation. This cancellation is revealed in Lemmas \ref{origen} and \ref{banda}.

Theorem \ref{teorema_gradiente} will be proved directly  as  a consequence  of the  decomposition of the resolvent of Navier-Lam\'e in terms of the resolvent of  the Laplacian  at two different  energies  given in \cite[Lemma 2.2]{BFPRV2012}. In the proof of \eqref{estimacion_gradu} we  work in the frequency domain.
We suspect  that  a  similar proof  could  be given  for  \eqref{estimacion_u}, but  so far  we  have not been able to do this.

\textbf{Notation.} 
From now on, all vectors will be considered column vectors, and all differential operators will be assumed to act on column vectors.


For non-negative
quantities $X$ and $Y$ we will use $X\lesssim Y$ ($X\gtrsim Y$) to
denote the existence of a positive constant $c$, depending on at
most $d$, such that $X\leq c\,Y$ ($X\geq c\,Y$). We will write
$X\sim Y$ if both $X\lesssim Y$ and $X\gtrsim Y$.

\section{The class of radial functions with bounded X-ray transform.}
\label{sec_contraejemplo}

In \cite{BFPRV2012}, the Rubio de Francia Algorithm was used to contruct  an $A_2$ majorant for functions in some  classes. This  was essential to transfer weighted  estimates from the Helmholt equation to the Navier-Lam\'e equation. It requires an a priori bound for the Hardy-Littlewood Maximal operator in the  class under consideration.   In this  section  we will see that the class of radial  functions with bounded  X-ray  transform behaves  badly under    the classical centered Hardy-Littlewood maximal operator. More precisely, we will prove the following proposition.
 
\begin{proposition}
\label{Contraejemplo}
The classical centered Hardy-Littlewood  maximal operator on $\mathbb R^d$ is not a bounded operator on  the radial Mizohata-Takeuchi class $\mathcal{T}_{rad}(\mathbb{R}^d)$.
\end{proposition}
\begin{proof}
 If we denote the classical centered Hardy-Littlewood  maximal operator by $M$, we will show that the following estimate does not hold:
\begin{equation}
    \label{falsa}
    |||Mf|||\lesssim |||f|||,
    \qquad\qquad
    f\in \mathcal{T}_{rad}(\mathbb{R}^d).
\end{equation}

If $f\in \mathcal{T}_{rad}(\mathbb{R}^d)$, then $f$ is a radial function. Writing $f(x)=f_0(|x|)$ and using polar coordinates, we have that
\begin{equation*}
 \begin{split}
    Mf(x)&=\sup_{s>0}\frac{1}{|B_s(x)|} \int_{B_s(x)}|f(y)|\,dy
    \\
    &\ge\sup_{0<s<|x|}\frac{1}{|B_s(x)|} \int_{B_s(x)}|f(y)|\,dy
    \\&\sim\sup_{0<s<|x|} \frac{1}{s^d}\int_{|x|-s}^{|x|+s}|f_0(t)|\,h(s,t)dt,
\end{split}   
\end{equation*}
where $B_s(x)$ is the ball of radius $s$ centered at $x$ and
\begin{equation*}
 h(s,t)=\int_{B_s(x)\cap \{|y|=t\}}d\sigma(\omega).
\end{equation*}
Taking into account that if $||x|-t|<s/2$, then $h(s,t)\sim s^{d-1}$, we have 
\begin{equation}
\label{Nf0}
    Mf(x)\gtrsim\mathcal{N}f_0(|x|):= \sup_{0< s<|x|}\frac{1}{2s}\int_{|x|-s}^{|x|+s}|f_0(t)|\,dt.
\end{equation}
From here, using \eqref{normaTakeuchi}, we see that
\begin{equation}
\label{norm_maximal}
    |||Mf|||\gtrsim 
    \sup_{\mu\ge0}\int_{\mu}^{\infty}\frac{r\,\mathcal{N}f_0(r)}{(r^2-\mu^2)^{1/2}}\,dr.
\end{equation}

Consider two small positive parameters, 
$\eta$ and $\delta$, satisfying $0<\eta^2\ll\delta\ll\eta\ll1$ and let
$N$ be a natural number satisfying $N\sim \eta^{-1}$. 
We introduce the following function depending on these parameters:
\begin{equation*}
    f(x)=f_0(|x|)=\sum_{j=0}^N\chi_{I_j}(|x|),
\end{equation*}
where $I_j=[2+j\eta  ,2+j\eta +\delta]$ for $1\le j\le N$.

On one hand, for this function $f$, we have that
\begin{equation*}
   |||f|||
   =\sup_{\mu\ge0}\int_{\mu}^{\infty}\frac{r\,f_0(r)}{(r^2-\mu^2)^{1/2}}\,dr
   \sim\sup_{2\le \mu\le 3}\int_{\mu}^{3}\frac{r\,f_0(r)}{(r^2-\mu^2)^{1/2}}\,dr. 
\end{equation*}
Fix $2\le \mu\le 3$. Then there exists $j_0\in\{0,\ldots,N-1\}$ such that $2+j_0\eta\le \mu< 2+(j_0+1)\eta$. And thus, we can write
\begin{equation}
\label{norm_f_contraejemplo}
   \int_{\mu}^{\infty}\frac{r\,f_0(r)}{(r^2-\mu^2)^{1/2}}\,dr \ \lesssim \
   \sum_{j=j_0}^N {\mathcal I}_j,
\end{equation}
where
\begin{equation}
\label{Ij0_contraejemplo}
    {\mathcal{I}}_{j_0}:=\int_{\mu}^{2+j_0\eta+\delta}\frac{r}{(r^2-\mu^2)^{1/2}}\,dr=\sqrt{(2+j_0\eta+\delta)^2-\mu^2}\lesssim \delta^{1/2},
\end{equation}
and for $j=j_0+1,\ldots,N,$
\begin{align}
    {\mathcal I}_j
    &:=\int_{2+j\eta}^{2+j\eta+\delta}\frac{r}{(r^2-\mu^2)^{1/2}}\,dr
    \nonumber
    \\
    &=\sqrt{(2+j\eta+\delta)^2-\mu^2}-\sqrt{(2+j\eta)^2-\mu^2}
    \nonumber
    \\
    &\lesssim \frac{\delta}{(j-j_0)^{1/2}\eta^{1/2}}.
    \label{Ij_contraejemplo}
\end{align}
Inserting \eqref{Ij0_contraejemplo} and \eqref{Ij_contraejemplo} in \eqref{norm_f_contraejemplo} and taking the supremum, we obtain
\begin{equation}
\label{norma_f0}
    |||f|||
    \lesssim 
    \delta^{1/2} + \frac{\delta}{\eta^{1/2}}\sum_{j=j_0+1}^N \frac{1}{(j-j_0)^{1/2}}
    \sim
    \delta\eta^{-1}
\end{equation}
since $\eta^2\ll\delta$.
On the other hand, from \eqref{norm_maximal}, we see that
\begin{equation}
\label{norm_maximal_f0}
    |||Mf|||
    \gtrsim 
    \int_2^3\frac{r\,\mathcal{N}f_0(r)}{(r^2-4)^{1/2}}\,dr
    \gtrsim
    \sum_{j=1}^N (j\eta)^{-1/2}
    \int_{2+j\eta}^{2+(j+1)\eta}\mathcal{N}f_0(r)\,dr.
\end{equation}
For a fixed $r\in [2+j \eta, 2+(j+1)\eta]$ with $j=1,\ldots,N$, we see from \eqref{Nf0} that
\begin{equation*}
    \mathcal{N}f_0(r)
    \ge \frac{1}{2(r-2-j\eta+\delta)}\int_{2+j\eta-\delta}^{2r-2-j\eta+\delta}|f_0(t)|\,dt
    \gtrsim
    \frac{\delta}{r-2-j\eta+\delta}.
\end{equation*}
Inserting this in \eqref{norm_maximal_f0} we obtain that
\begin{align}
    \label{norm_maximal_f0_fin}
    |||Mf|||
    &\gtrsim 
    \frac{\delta}{\eta^{1/2}}\sum_{j=1}^N \frac{1}{{\sqrt{j}}}
    \int_{2+j\eta}^{2+(j+1)\eta}\frac{dr}{r-2-j\eta+\delta}
    \nonumber
    \\
    &=\frac{\delta}{\eta^{1/2}}\ln{\left(\frac{\eta+\delta}{\delta}\right)}\sum_{j=1}^N  \frac{1}{{\sqrt{j}}}
    \nonumber
    \\
    &\sim
    \delta \eta^{-1}\ln{\left(\frac{\eta}{\delta}\right)}.
\end{align}
 
Hence, if \eqref{falsa} held, from \eqref{norma_f0} and \eqref{norm_maximal_f0_fin}, we would conclude that
\begin{equation*}
  \ln{\left(\frac{\eta}{\delta}\right)}\lesssim 1,  
\end{equation*}
which is impossible since $0<\delta\ll\eta\ll1$.
\end{proof}

\section{The fundamental solution in spherical harmonics.}
The fundamental solution of equation \eqref{ecuacion} satisfying the outgoing Kupradze conditions \eqref{radiacionup} and \eqref{radiacionus} is an $d\times d$ matrix. 
In this section we present series expansions of this fundamental solution in dimensions 2 and 3. 
These expansions will be key in proving  Theorem \ref{principal}.
We will first focus on dimension 2, and then we will address dimension 3, which is much more complicated.

In dimension 2, the fundamental solution can be written as (see \cite{Arens})
\begin{equation}
    \label{fundamental_solution}
    \Phi(x,y)
    =\frac{i}{4 \omega^2}\left(\nabla_x\nabla_y^t  H_0^{(1)}(k_p|x-y|)
    +\nabla_x^\perp\left(\nabla_y^\perp\right)^t H_0^{(1)}(k_s|x-y|)\right),
\end{equation}
where $x=(x_1,x_2)^t,$ $\nabla_x=(\partial_{x_1},\partial_{x_2})^t$, $\nabla_x^\perp=(-\partial_{x_2},\partial_{x_1})^t$ and $H_0^{(1)}$ denotes the Hankel function of the first kind and order zero.

Writing in polar coordinates $x=(r \cos \theta , r \sin \theta)^t$ and $y=(t \cos \varphi , t \sin \varphi)^t$, we have the following addition formula (see \cite[p.75]{ColtonKress})
\begin{equation}\label{addition}
H_0^{(1)}(k| x-y |)= \sum_{n \in \mathbb{Z}} H_{n}^{(1)}
(kr) J_{n}(kt) e^{-in \theta}e^{in \varphi}, \qquad k>0, \ |x|>|y|,
\end{equation}
where $J_{n}$ and $H_{n}^{(1)}$ are the Bessel functions of the first and third kind, respectively. 

The series and its term-by-term first derivatives with respect to $|x|$ and $|y|$ are absolutely and uniformly convergent on compact subsets of $|x|>|y|$.

The expression of the fundamental solution given in \eqref{fundamental_solution} along with the addition formula \eqref{addition} will allow us to obtain the following addition formula for the fundamental matrix solution.

\begin{proposition}
\label{propo_addition_2D}
If we write in polar coordinates $x=(r\cos\theta,r\sin\theta)$ and $y=(t\cos\varphi,t\sin\varphi)$ then, for $|x|>|y|$ the fundamental solution of \eqref{ecuacion} in dimension 2 is a $2\times2$ matrix that can be written as
\begin{equation}
\label{addition_expansion_2d}
    \Phi(x,y)=\Phi^1(x,y)+\Phi^2(x,y)+\Phi^3(x,y),
\end{equation}
with
\begin{align}
    \Phi^1(x,y)&=
    \frac{i}{8\omega^2}\sum _{n\in \mathbb Z}
    H_{n,n}^+(r,t)\mathrm{I}\,e^{in\theta}e^{-in\varphi},
    \label{Phi1}
    \\
    \Phi^2(x,y)&=
    \frac{1}{16\omega^2}\sum _{n\in \mathbb Z}
    H_{n-2,n}^-(r,t)\mathrm{A}\,e^{i(n-2)\theta}e^{-in\varphi},
\label{Phi2}
    \\
    \Phi^3(x,y)&=
    \frac{1}{16\omega^2}\sum _{n\in \mathbb Z}
    H_{n+2,n}^-(r,t)\mathrm{B}\,e^{i(n+2)\theta}e^{-in\varphi},
\label{Phi3}
\end{align}
where for $n\in\mathbb{Z}$ and $m=n-2$ or $m=n+2$,
\begin{align}
    \label{H+}
    H_{n,n}^+( r,t)&:=k_p^2H_n^{(1)}(k_pr)J_n(k_pt)+k_s^2H_n^{(1)}(k_sr)J_n(k_st),
    \\
    \label{H-}
    H_{m,n}^-(r,t)&:=
    k_p^2H_m^{(1)}(k_pr)J_n(k_pt)-k_s^2H_m^{(1)}(k_sr)J_n(k_st),
\end{align}
$\mathrm{I}$ is the identity matrix and
\begin{equation}
    \label{matrices2D}
    A=
    \begin{pmatrix}
    -i & -1 \\ -1 & i
    \end{pmatrix},
    \qquad
    B=
    \begin{pmatrix}
    -i & 1\\ 1 & i
    \end{pmatrix}.
\end{equation}
The series and its term-by-term first derivatives with respect to $|x|$ and
$|y|$ are absolutely and uniformly convergent on compact subsets of $|x| > |y|$.
\end{proposition}

For convenience, from now on, we will write $H_{n}$ instead of $H_{n}^{(1)}$.

\emph{Proof of Proposition \ref{propo_addition_2D}}.
Following \cite{Arens}, we define the Navier eigenfunctions
\begin{equation*}
    \mathbf{e}_{p,n}(x):=\nabla_x(J_n(k_pr)e^{in\theta}),\qquad
    \mathbf{e}_{s,n}(x):=\nabla_x^\perp(J_n(k_sr)e^{in\theta}),
\end{equation*}
with $n\in\mathbb{Z}, x\in\mathbb{R}^2$ and, the radiating Navier eigenfunctions
\begin{equation*}
    \mathbf{r}_{p,n}(x):=\nabla_x(H_n(k_pr)e^{-in\theta}),\qquad
    \mathbf{r}_{s,n}(x):=\nabla_x^\perp(H_n(k_sr)e^{-in\theta}),
\end{equation*}
with $n\in\mathbb{Z}, x\in\mathbb{R}^2\setminus\{0\}.$

With these definitions, from \eqref{fundamental_solution} and \eqref{addition}, we get the following addition formula for the fundamental solution of the Lamé equation (see \cite{Arens}): if $|x|>|y|$, then
\begin{equation}
    \label{addition_Lame_2D}
    \Phi(x,y)
    =\frac{i}{4 \omega^2}
    \sum_{n \in \mathbb{Z}}
    (\mathrm{r}_{p,n}(x)\mathrm{e}_{p,n}^t(y)+
    \mathrm{r}_{s,n}(x)\mathrm{e}_{s,n}^t(y)),
\end{equation}
where the series converges uniformly on compact subsets.

For convenience, we define 
\begin{equation}
\label{JknHkn}
    J_{k,n}(x):=J_n(kr)e^{in\theta},
    \ \  {\rm and} \ \
    H_{k,n}(x):=H_n(kr)e^{-in\theta}.
\end{equation}
And we write
$\Phi=(\phi_{\ell j})$ with
\begin{equation}
\label{phi_componente}
    \phi_{\ell j}=\frac{i}{4 \omega^2}\sum_{n \in \mathbb{Z}}\phi_{\ell j}^n.
\end{equation}
Then, from \eqref{addition_Lame_2D} we have 
\begin{align}
\label{phi11_componente}
    \phi_{11}^n(x,y)&=
    \partial_{x_1}H_{k_p,n}(x)\partial_{y_1}J_{k_p,n}(y)
    +\partial_{x_2}H_{k_s,n}(x)\partial_{y_2}J_{k_s,n}(y),
    \\
\label{phi22_componente}
    \phi_{22}^n(x,y)&=
    \partial_{x_2}H_{k_p,n}(x)\partial_{y_2}J_{k_p,n}(y)
    +\partial_{x_1}H_{k_s,n}(x)\partial_{y_1}J_{k_s,n}(y),
    \\
\label{phi12_componente}
    \phi_{\ell j}^n(x,y)&=
    \partial_{x_\ell}H_{k_p,n}(x)\partial_{y_j}J_{k_p,n}(y)
    -\partial_{x_j}H_{k_s,n}(x)\partial_{y_\ell}J_{k_s,n}(y),
    \qquad \ell\neq j.
\end{align}

On the other hand, using the chain rule, for $j=1,2$ we have 
\begin{align}
    \partial_{x_j}J_{k,n}(x)&=kJ'_n(kr)e^{in\theta}\partial_{x_j}r+inJ_n(kr)e^{in\theta}\partial_{x_j}\theta,
\label{parcial_Jkn}
    \\
    \partial_{x_j}H_{k,n}(x)&=kH_n'(kr)e^{-in\theta}\partial_{x_j}r-inH_n(kr)e^{-in\theta}\partial_{x_j}\theta.
\label{parcial_Hkn}
\end{align}

Taking into account that
\begin{align*}
    \partial_{x_1}r&=\frac{e^{i\theta}+e^{-i\theta}}{2},
    &\partial_{x_1}\theta&=\frac{i(e^{i\theta}-e^{-i\theta})}{2r},
    \\
    \partial_{x_2}r&=\frac{-i(e^{i\theta}-e^{-i\theta})}{2},
    &\partial_{x_2}\theta&=\frac{e^{i\theta}+e^{-i\theta}}{2r},
\end{align*}
from \eqref{parcial_Jkn}, we get
\begin{align}
    \partial_{x_1}J_{k,n}(x)=&\frac{k}{2}\left[J'_n(kr)(e^{i(n+1)\theta}+e^{i(n-1)\theta})
    -\frac{nJ_n(kr)}{kr}(e^{i(n+1)\theta}-e^{i(n-1)\theta})\right],
\label{parcial1_Jkn}
    \\
    \partial_{x_2}J_{k,n}(x)=&-i\frac{k}{2}\left[J'_n(kr)(e^{i(n+1)\theta}-e^{i(n-1)\theta})
    -\frac{nJ_n(kr)}{kr}(e^{i(n+1)\theta}+e^{i(n-1)\theta})\right],
\label{parcial2_Jkn}
\end{align}
and from \eqref{parcial_Hkn}, we obtain
\begin{align}
    \partial_{x_1}H_{k,n}(x)=\frac{k}{2}&\left[H'_n(kr)(e^{-i(n-1)\theta}+e^{-i(n+1)\theta})\right.
\nonumber
    \\
    &+\left.\frac{nH_n(kr)}{kr}(e^{-i(n-1)\theta}-e^{-i(n+1)\theta})\right],
\label{parcial1_Hkn}
    \\
    \partial_{x_2}H_{k,n}(x)=-i\frac{k}{2}&\left[H'_n(kr)(e^{-i(n-1)\theta}-e^{-i(n+1)\theta})\right.
\nonumber
    \\
    &+\left.\frac{nH_n(kr)}{kr}(e^{-i(n-1)\theta}+e^{-i(n+1)\theta})\right].
\label{parcial2_Hkn}
\end{align}

If $G_n(r)$ is $J_n(r)$ or $H_n(r)$, the following recurrence relations hold for all $n\in\mathbb{Z}$ (see \cite[pp.361, identities in 9.1.27]{AS1967}:
\begin{align*}
    2G'_n(r)&=G_{n-1}(r)-G_{n+1}(r),\\
    \frac{2n}{r}G_n(r)&=G_{n-1}(r)+G_{n+1}(r).
\end{align*}

Using these recurrence relations in \eqref{parcial1_Jkn}-\eqref{parcial2_Hkn}, and grouping terms we obtain
\begin{align}
\label{parcial1J}
    \partial_{x_1}J_{k,n}(x)=&\frac{k}{2}\left(J_{k,n-1}(x)-J_{k,n+1}(x)\right),
    \\
\label{parcial2J}
    \partial_{x_2}J_{k,n}(x)=&i\frac{k}{2}\left(J_{k,n-1}(x)+J_{k,n+1}(x)\right),
\\
\label{parcial1H}
    \partial_{x_1}H_{k,n}(x)=&\frac{k}{2}\left(H_{k,n-1}(x)-H_{k,n+1}(x)\right), \ \ {\rm and}
    \\
\label{parcial2H}
    \partial_{x_2}H_{k,n}(x)=&-i\frac{k}{2}\left(H_{k,n-1}(x)+H_{k,n+1}(x)\right).
\end{align}

We start by calculating the term $\phi_{11}$ of the fundamental matrix solution $\Phi$.
From \eqref{phi_componente} with $\ell=j=1$ we see that
\begin{equation}
\label{phi11}
    \phi_{11}=\frac{i}{4 \omega^2}\sum_{n \in \mathbb{Z}}\phi_{11}^n,
\end{equation}
where $\phi_{11}^n$ is given in \eqref{phi11_componente}. 

Inserting the identities \eqref{parcial1J}-\eqref{parcial2H} in \eqref{phi11_componente} we get
\begin{align}
    \phi_{11}^n(x,y)&=
    \frac{k_p^2}{4}\left[
    H_{k_p,n-1}(x)J_{k_p,n-1}(y)
    +H_{k_p,n+1}(x)J_{k_p,n+1}(y)
    \right]
\nonumber
    \\
    &-
    \frac{k_p^2}{4}\left[
    H_{k_p,n-1}(x)J_{k_p,n+1}(y)
    +H_{k_p,n+1}(x)J_{k_p,n-1}(y)
    \right]
\nonumber
    \\
    &+
    \frac{k_s^2}{4}\left[
    H_{k_s,n-1}(x)J_{k_s,n-1}(y)
    +H_{k_s,n+1}(x)J_{k_s,n+1}(y)
    \right]
\nonumber
    \\
    &+
    \frac{k_s^2}{4}\left[
    H_{k_s,n-1}(x)J_{k_s,n+1}(y)
    +H_{k_s,n+1}(x)J_{k_s,n-1}(y)
    \right].
\label{phi11_n}
\end{align}

Now we make the change of variables in the summation index $n$ appearing in \eqref{phi11}, that are necessary for the terms $J_{k,m}(y)$ with $k=k_p$ or $k=k_s$ and $m=n-1$ or $m=n+1$ to become $J_{k,n}(y)$.
After these changes of variables, grouping terms, we obtain
\begin{align*}
   \phi_{11}(x,y)&=
   \frac{i}{8 \omega^2}\sum_{n \in \mathbb{Z}}
   \left[
    k_p^2
    H_{k_p,n}(x)J_{k_p,n}(y)
    +k_s^2
    H_{k_s,n}(x)J_{k_s,n}(y)
    \right]
\nonumber
    \\
    &-
    \frac{i}{16 \omega^2}\sum_{n \in \mathbb{Z}}
    k_p^2\left[
    H_{k_p,n-2}(x)J_{k_p,n}(y)
    +H_{k_p,n+2}(x)J_{k_p,n}(y)
    \right]
\nonumber
    \\
    &+
    \frac{i}{16 \omega^2}\sum_{n \in \mathbb{Z}}
    k_s^2\left[
    H_{k_s,n-2}(x)J_{k_s,n}(y)
    +H_{k_s,n+2}(x)J_{k_s,n}(y)
    \right]. 
\end{align*}

Using \eqref{JknHkn}, \eqref{H+} and \eqref{H-} in this expression,
we get
\begin{align*}
   \phi_{11}&=
   \frac{i}{8 \omega^2}\sum_{n \in \mathbb{Z}}
    H_{n,n}^+(r,t)
    e^{-in\theta}e^{in\varphi}
\nonumber
    \\
    &-
    \frac{i}{16 \omega^2}\sum_{n \in \mathbb{Z}}
    H_{n-2,n}^-(r,t)
    e^{-i(n-2)\theta}e^{in\varphi}
\nonumber
    \\
    &-
    \frac{i}{16 \omega^2}\sum_{n \in \mathbb{Z}}
    H_{n+2,n}^-(r,t)
    e^{-i(n+2)\theta}e^{in\varphi}.
\end{align*}
Or equivalently,
\begin{align}
   \phi_{11}&=
   \frac{i}{8 \omega^2}\sum_{n \in \mathbb{Z}}
    H_{n,n}^+(r,t)
    e^{in\theta}e^{-in\varphi}
\nonumber
    \\
    &-
    \frac{i}{16 \omega^2}\sum_{n \in \mathbb{Z}}
    H_{n-2,n}^-(r,t)
    e^{i(n-2)\theta}e^{-in\varphi}
\nonumber
    \\
    &-
    \frac{i}{16 \omega^2}\sum_{n \in \mathbb{Z}}
    H_{n+2,n}^-(r,t)
    e^{i(n+2)\theta}e^{-in\varphi}.
\label{phi11b}
\end{align}
To obtain this equivalence, we used the fact that for all $n\in\mathbb{Z}$, 
\begin{equation}
\label{cambio_signo}
    H_{-n-2}J_{-n}=H_{n+2}J_n,
    \qquad H_{-n}J_{-n}=H_{n}J_n, 
    \qquad H_{-n+2}J_{-n}=H_{n-2}J_n.
\end{equation}

Arguing in a similar way one can see that
\begin{align}
   \phi_{22}&=
   \frac{i}{8 \omega^2}\sum_{n \in \mathbb{Z}}
    H_{n,n}^+(r,t)
    e^{in\theta}e^{-in\varphi}
\nonumber
    \\
    &+
    \frac{i}{16 \omega^2}\sum_{n \in \mathbb{Z}}
    H_{n-2,n}^-(r,t)
    e^{i(n-2)\theta}e^{-in\varphi}
    \nonumber
    \\
    &+
    \frac{i}{16 \omega^2}\sum_{n \in \mathbb{Z}}
    H_{n+2,n}^-(r,t)
    e^{i(n+2)\theta}e^{-in\varphi},
\label{phi22}
\end{align} 
and that 
\begin{align}
   \phi_{12}&=
    \frac{1}{16 \omega^2}\sum_{n \in \mathbb{Z}}
    H_{n-2,n}^-(r,t)
    e^{i(n-2)\theta}e^{-in\varphi}
    \nonumber
    \\
    &-
   \frac{1}{16 \omega^2}\sum_{n \in \mathbb{Z}}
    H_{n+2,n}^-(r,t)
    e^{i(n+2)\theta}e^{-in\varphi}.
\label{phi12}
\end{align}
On the other hand, from \eqref{fundamental_solution} we have that $\Phi$ is a symmetric matrix, and therefore
\begin{equation}
    \label{phi21}
    \phi_{21}=\phi_{12}.
\end{equation}

Identity \eqref{addition_expansion_2d} follows from \eqref{phi11b}-\eqref{phi21}.
\hfill $\square$

\begin{remark}
    \label{signo_exponencial}
    Notice that due to \eqref{cambio_signo}, the series expansion \eqref{addition_expansion_2d} is also true if we substitute $e^{in\theta}$ and $e^{-in\varphi}$ for $e^{-in\theta}$ and $e^{in\varphi}$, respectively.
\end{remark}

In dimension 3, the series expansion analogous to \eqref{addition_expansion_2d} is much more complicated. In this case, we can write the fundamental solution as (see \cite{Arens})
\begin{equation}
\label{fundamental3D}
\mathrm{\Phi}(x,y)=\frac{1}{4\pi\omega^2}
\left(\nabla_x\nabla_y^t \frac{e^{ik_p|x-y|}}{|x-y|}
+
\mathcal{D}_x\mathcal{D}_y^t
\frac{e^{ik_s|x-y|}}{|x-y|}
\right),
\qquad x\neq y,
\end{equation}
where $x=(x_1,x_2,x_3)^t,$ $\nabla_x=(\partial_{x_1},\partial_{x_2},\partial_{x_3})^t$ and
\begin{equation*}
\mathcal{D}:=\left(
\begin{array}{ccc}
0&-\partial_{x_3}&\partial_{x_2}\\
\partial_{x_3}&0&-\partial_{x_1}\\
-\partial_{x_2}&\partial_{x_1}&0\\
\end{array}
\right).
\end{equation*}

The series expansion in the 3-dimensional case involves spherical harmonics.
We consider the following normalization for them:
\begin{equation}
\label{armonico}
Y_n^m:=\gamma_n^me^{im\varphi}P_n^{m}(\cos\theta),
\qquad
-n\le m\le n,\ n=0,1,\ldots,
\end{equation}
where
\begin{equation}
\label{gamma_n^m}
\gamma_n^m:=\sqrt{\frac{(2n+1)(n-m)!}{4\pi(n+m)!}}.
\end{equation}
Here $P_n^{m}$ is the associated Legendre function of degree $n$ and order $m$ defined by
\begin{equation}
\label{legendre}
P_n^{m}(z):=(-1)^m(1-z^2)^\frac{m}{2}\frac{d^mP_n(z)}{dz^m}, 
\qquad
z\in[-1,1],\ 0\le m\le n,
\end{equation}
where $P_n$ is the Legendre polynomial of degree $n$.
And we use the standard convention that for $m>0$,
\begin{equation}
\label{Pnm_negativo}
P_n^{-m}(z):=(-1)^m\frac{(n-m)!}{(n+m)!}\,P_n^m(z).
\end{equation}

The following result can be found in \cite{Addition_Lame_3D}.
\begin{theorem}
    \label{adition_theorem}
    Let $Y_n^m$, $m=-n,\ldots,n$, $n=0,1,\ldots$, be the set of orthonormal spherical harmonics given in \eqref{armonico}. 
    And denote the spherical coordinates of $x$ and $y$ by $(r,\theta_x,\varphi_x)$ and $(t,\theta_y,\varphi_y)$ respectively.
	Then for $|x|>|y|$, the fundamental solution of \eqref{ecuacion} is a matrix that can be written as
 \begin{equation}
 \label{additionLame}
    \Phi(x,y)=\mathrm{\Psi}(x,y)+\mathrm{\Phi}_-(x,y)+\mathrm{\Phi}_0(x,y)+\mathrm{\Phi}_+(x,y),
\end{equation}
with
    \begin{align}
	\mathrm{\Psi}(x,y)&=
    \frac{1}{\omega^2}\sum_{n=0}^\infty
    \left(h_{n,n}^0(k_s,r,t)\mathrm{S}+h_{n,n}^0(k_p,r,t)\mathrm{P}\right)
    \sum_{m=-n}^{n}Y_{n}^{m}(\theta_x,\varphi_x)\overline{Y_n^m}(\theta_y,\varphi_y),
\label{Phi1_3D}
    \\
    \mathrm{\Phi}_-(x,y)&=
    \frac{1}{\omega^2}\sum_{n=0}^\infty
    h_{n-2,n}^-(r,t)
    \sum_{m=-n}^{n}
    S_{-,n}^m(\theta_x,\varphi_x)\overline{Y_n^m}(\theta_y,\varphi_y),
\label{Phi2_3D}
    \\
    \mathrm{\Phi}_0(x,y)&=
    \frac{1}{\omega^2}\sum_{n=0}^\infty
    h_{n,n}^-(r,t)
    \sum_{m=-n}^{n}
     S_{0,n}^m(\theta_x,\varphi_x)\overline{Y_n^m}(\theta_y,\varphi_y),
\label{Phi3_3D}
    \\
    \mathrm{\Phi}_+(x,y)&=
    \frac{1}{\omega^2}\sum_{n=0}^\infty
    h_{n+2,n}^-(r,t)
    \sum_{m=-n}^{n}
     S_{+,n}^m(\theta_x,\varphi_x)\overline{Y_n^m}(\theta_y,\varphi_y),
\label{Phi4_3D}
    \end{align}
where 
\begin{align}
\label{h0}
h_{n,n}^0(k,r,t)&:=
k^3h_{n}^{(1)}(kr)j_n(kt),
\\
\label{h-}
h_{n_1,n}^-(r,t)&:=
k_p^3h_{n_1}^{(1)}(k_pr)j_n(k_pt)
-k_s^3h_{n_1}^{(1)}(k_sr)j_n(k_st),
\end{align}
with $j_n$ and $h_n^{(1)}$ denoting the spherical Bessel functions of the first and third kind respectively,
\begin{equation}
    \label{SP}
    \mathrm{S}=\left(
    \begin{array}{ccc}
       1&0&0  \\
       0&1&0  \\
       0&0&0 
    \end{array}
    \right),
    \qquad
    \mathrm{P}=\left(
    \begin{array}{ccc}
       0&0&0  \\
       0&0&0  \\
       0&0&1 
    \end{array}
    \right),
\end{equation}
\begin{align}
    S_{-,n}^m(\theta_x,\varphi_x)
    =&
    \mathrm{A}_{-,n}^{m}
    Y_{n-2}^{m-2}(\theta_x,\varphi_x)
    +
    \mathrm{B}_{-,n}^{m}
    Y_{n-2}^{m-1}(\theta_x,\varphi_x)
    +
    \mathrm{C}_{-,n}^{m}
    Y_{n-2}^{m}(\theta_x,\varphi_x)
    \nonumber
    \\
    &+
    \mathrm{D}_{-,n}^{m}
    Y_{n-2}^{m+1}(\theta_x,\varphi_x)
    +
    \mathrm{E}_{-,n}^{m}
    Y_{n-2}^{m+2}(\theta_x,\varphi_x),   
\end{align}
\begin{align}
    S_{0,n}^m(\theta_x,\varphi_x)
    =&
    \mathrm{A}_{0,n}^{m}
    Y_{n}^{m-2}(\theta_x,\varphi_x)
    +
    \mathrm{B}_{0,n}^{m}
    Y_{n}^{m-1}(\theta_x,\varphi_x)
    +
    \mathrm{C}_{0,n}^{m}
    Y_{n}^{m}(\theta_x,\varphi_x)
    \nonumber
    \\
    &+
    \mathrm{D}_{0,n}^{m}
    Y_{n}^{m+1}(\theta_x,\varphi_x)
    +
    \mathrm{E}_{0,n}^{m}
    Y_{n}^{m+2}(\theta_x,\varphi_x),   
\end{align}
\begin{align}
    S_{+,n}^m(\theta_x,\varphi_x)
    =&
    \mathrm{A}_{+,n}^{m}
    Y_{n+2}^{m-2}(\theta_x,\varphi_x)
    +
    \mathrm{B}_{+,n}^{m}
    Y_{n+2}^{m-1}(\theta_x,\varphi_x)
    +
    \mathrm{C}_{+,n}^{m}
    Y_{n+2}^{m}(\theta_x,\varphi_x)
    \nonumber
    \\
    &+
    \mathrm{D}_{+,n}^{m}
    Y_{n+2}^{m+1}(\theta_x,\varphi_x)
    +
    \mathrm{E}_{+,n}^{m}
    Y_{n+2}^{m+2}(\theta_x,\varphi_x),   
\label{S+}
\end{align}
with $\mathrm{A}_{\#,n}^{m},\mathrm{B}_{\#,n}^{m},\mathrm{C}_{\#,n}^{m},\mathrm{D}_{\#,n}^{m}$ and $\mathrm{E}_{\#,n}^{m}$ for $\#=-,+,0,$ being $3\times 3$ matrices of the discrete variables $n$ and $m$, with $n\ge0,$ and $-n\le m\le n$ uniformly bounded.

The series and its term-by-term first derivatives with respect to $|x|$ and $|y|$ are absolutely and uniformly convergent on compact subsets of $|x| > |y|$.
\end{theorem}
\begin{remark}
\label{r<t}
    It follows from \eqref{fundamental_solution} that 
    $\Phi(x,y)=\Phi(y,x)$.
    Hence in the region $|x| < |y|$,  \eqref{addition_expansion_2d} holds by changing the roles of $(r,\theta)$ by $(t,\varphi)$.

    Analogously, from \eqref{fundamental3D} we have that in the region $|x| < |y|$,  \eqref{additionLame} holds by changing the roles of $(r,\theta_x,\varphi_x)$ by $(t,\theta_y,\varphi_y)$.
\end{remark}
\section{Proof of Theorem \ref{principal}}
Before proving the theorem, we will state several lemmas that will be crucial in its proof. These lemmas provide estimates for certain integrals where the integrands involve Bessel and Hankel functions. Their proofs are postponed until Section \ref{demos_lemas}.

First we present a lemma proved by Barceló, Ruiz and Vega in \cite{BRV1997}.
\begin{lemma}\cite[Lemma 6]{BRV1997}
    \label{lema_BRV}
    Let $\mu\ge 1/2$ and $d\ge 2$. If $V\in \mathcal{T}_{rad}(\mathbb{R}^d)$, then
    \begin{equation}
        \label{est_BRV}
        \int_0^\infty\left|J_\mu(r)\right|^2V(r)r\,dr
        +
        \int_{\mu-\mu^{\frac{1}{3}}}^\infty\left|H^{(1)}_\mu(r)\right|^2V(r)r\,dr
        \lesssim |||V|||.
    \end{equation}
\end{lemma}
In the following lemma we give an estimate where the constant that appears depends on the orders of the Bessel functions. 
\begin{lemma}
    \label{mu_peque}
    Let $\mu=n$ or $\mu=n+\frac{1}{2}$ with $n$ a non-negative integer, and let $\nu=\mu-2$ or $\nu=\mu+2$, $a>1$ and $d\ge 2$. If $V\in \mathcal{T}_{rad}(\mathbb{R}^d)$, then
    \begin{equation}
\label{est_mu_peque}
    \int_0^\infty \int_r^\infty \left|H_\nu^{(1)}(t)J_\mu(r)-aH_\nu^{(1)}(at)J_\mu(ar)\right|^2 tV(t) dt\,rV(r) dr
    \le c_{\mu,a} |||V|||^2,
\end{equation}
where $c_{\mu,a}$ is a constant depending on $\mu$ and $a$.
\end{lemma}
The rest of the lemmas we present give estimates when the orders of the Bessel functions are sufficiently large.
\begin{lemma}
    \label{origen}
    Let $a> 1$ and $d\ge 2$. For $\mu>\mu_0$ with $\mu_0=\mu_0(a)$ large enough,  
    if $V\in \mathcal{T}_{rad}(\mathbb{R}^d)$, then
    \begin{equation}
        \label{est_origen}
        \int_0^{1+\frac{1}{2a}}\int_r^{1+\frac{1}{2a}}
        \left|J_\mu(r)H^{(1)}_{\mu+2}(t)-a^2J_\mu(ar)H^{(1)}_{\mu+2}(at)\right|^2
        V(t)t\,dt\,V(r)r\,dr
        \lesssim |||V|||^2.
    \end{equation}
\end{lemma}
\begin{lemma}
    \label{banda}
    Let $a> 1$ and $d\ge 2$. For $\mu>\mu_0$ with $\mu_0=\mu_0(a)$ large enough, 
    if $V\in \mathcal{T}_{rad}(\mathbb{R}^d)$, then 
    \begin{equation}
        \label{est_banda}
        \int_{1-\frac{1}{2a}}^{\frac{\mu}{2a}}\int_r^{r+\frac{1}{a}}
        \left|J_\mu(r)H^{(1)}_{\mu+2}(t)-a^2J_\mu(ar)H^{(1)}_{\mu+2}(at)\right|^2
        V(t)t\,dt\,V(r)r\,dr
        \lesssim |||V|||^2.
    \end{equation}
\end{lemma}
\begin{lemma}
    \label{escalones}
    Let $a> 1$ and $d\ge 2$. For $\mu$ large enough, let $N=N(\mu,a)=[\mu-2a]$, where $[\mu-2a]$ denotes the integer part of $\mu-2a$.
    If $V\in \mathcal{T}_{rad}(\mathbb{R}^d)$, then
    \begin{equation*}
        \label{est_escalones}
        \sum_{m=0}^{N}
        \mathrm{L}_m
        \lesssim |||V|||^2,
    \end{equation*}
where for $m=0,\ldots,N-1,$
    \begin{eqnarray*}
        \label{int_m_escalones}
        \mathrm{L}_m
        &=&
       \int_0^{1+\frac{m-1}{2a}}\left|J_\mu(r)\right|^2V(r)r\,dr
        \int_{1+\frac{m}{2a}}^{1+\frac{m+1}{2a}}\left|H^{(1)}_{\mu+2}(t)\right|^2V(t)t\,dt,
        \\
        \label{int_{|n|}_escalones}
        \mathrm{L}_{N}
        &=&
        \int_0^{1+\frac{N-1}{2a}}\left|J_\mu(r)\right|^2V(r)r\,dr
        \int_{1+\frac{N}{2a}}^{\infty}\left|H^{(1)}_{\mu+2}(t)\right|^2V(t)t\,dt.
    \end{eqnarray*}
\end{lemma}
\begin{lemma}
    \label{Juan}
    Let $a> 1$ and $d\ge 2$.
    For $\mu>\mu_0$ with $\mu_0=\mu_0(a)$ large enough, 
    if $V\in \mathcal{T}_{rad}(\mathbb{R}^d)$, then
    \begin{equation}
        \label{est_Juan}
        \int_{\frac{\mu}{2a}-\frac{1}{a}}^{\mu+2-(\mu+2)^{\frac{1}{3}}}\int_r^{\mu+2-(\mu+2)^{\frac{1}{3}}}
        |H^{(1)}_{\mu+2}(t)|^2
        V(t)t\,dt\,|J_\mu(r)|^2V(r)r\,dr
        \lesssim |||V|||^2.
    \end{equation}
\end{lemma}
\begin{lemma}
    \label{bueno}
    Let $d\ge 2$. For $\mu$ large enough, 
    if $V\in \mathcal{T}_{rad}(\mathbb{R}^d)$, then
    \begin{equation}
        \label{est_buena}
        \int_0^{\mu-\mu^{\frac{1}{3}}}\int_r^{\mu-\mu^{\frac{1}{3}}}
        |H^{(1)}_\mu(t)|^2
        V(t)t\,dt\,|J_{\mu+2}(r)|^2V(r)r\,dr
        \lesssim |||V|||^2.
    \end{equation}
\end{lemma}

\begin{remark}
    \label{caso_facil}
We would like to point out that in Lemma \ref{bueno}, unlike in Lemma \ref{Juan}, the Bessel function is of higher order, $\mu+2$, while the Hankel function is of lower order, $\mu$. This allows the lower limit of the integral in $r$ to be $0$ in Lemma \ref{bueno}, whereas in Lemma \ref{Juan} it can only be lowered to something on the order of $\mu/2$ (see Remark \ref{explica_bueno}).
\end{remark}

\begin{remark}
    Notice that in the previous lemmas the dimension $d$ does not play any role in the case of radial potentials $V$. The reason for this is that condition \eqref{normaTakeuchi} is independent of the dimension $d$ in this case.
\end{remark}

\emph{Proof of Theorem \ref{principal}.}
Note that it is sufficient to prove the theorem for the radial majorant of $V$ given in \eqref{mayorante_radial}, which, abusing notation, we will also denote by $V$.

For $\omega>0$ we introduce 
   \begin{equation}
    \label{V_a}
    V_\omega(r)=V\left(\frac{r}{\omega}\right)
    \end{equation}
    so that $V_\omega\in\mathcal{T}_{rad}(\mathbb{R}^d)$ and
    \begin{equation}
    \label{normaV_a}
    |||V_\omega|||=\omega\,|||V|||.
    \end{equation}
This allows us to reduce the proof of the theorem to the case where $\omega=1$ due to homogeneity. 
In fact, suppose we have proved the estimate \eqref{estimacion_u} for $\omega=1$; let us now see how the case of arbitrary $\omega$ follows.
    Given  $\mathbf{f}\in C_0^\infty(\mathbb{R}^d)$ let $\mathbf{u}$ be the unique solution of \eqref{ecuacion} satisfying the outgoing Kupradze radiation conditions \eqref{radiacionup} and \eqref{radiacionus}. Then, the function $\mathbf{u_\omega}(x)=\mathbf{u}\left( \frac{x}{\omega}  \right)$ satisfies the equation
    \begin{equation*}
    \label{reduccion_w_1}
    \Delta^* \mathbf{u_\omega}+\mathbf{u_\omega}=\frac{1}{\omega^2}\,\mathbf{f}_\omega,
    \end{equation*}
    where $\mathbf{f}_\omega(x)=\mathbf{f} \left( \frac{x}{\omega}  \right)$.
    Moreover, 
    $(\mathbf{u_\omega})_p(x)=\mathbf{u }_p \left(  \frac{x}{\omega}  \right) $ and 
    $(\mathbf{u_\omega})_s(x)=\mathbf{u }_s \left(  \frac{x}{\omega}  \right) $ 
    satisfy the outgoing Kupradze radiation conditions given in \eqref{radiacionup} and \eqref{radiacionus} with $k_p=\frac{1}{\sqrt{2 \mu + \lambda}}$ and $k_s=\frac{1}{\sqrt{ \mu }}$, respectively.
    Therefore,
    \begin{equation*}
	\label{estimacion_u_omega}
	\int_{\mathbb{R}^d}|\mathbf{u}_\omega(x)|^2V_\omega(x)dx
	\lesssim\frac{|||V_\omega|||^2}{\omega^4}
	\int_{\mathbb{R}^d}|\mathbf{f}_\omega(x)|^2V_\omega^{-1}(x)dx.
	\end{equation*}
    From here, we get \eqref{estimacion_u} for arbitrary $\omega$ making a change of variable and using \eqref{normaV_a}.

For simplicity, we will only write the proof in the 2 dimensional case and then explain the differences that appear in the 3 dimensional case.

\underline{2-D CASE}.
Consider $\mathbf{f}\in \mathcal{C}_0^{\infty}(\mathbb{R}^2)$. Writing $y=(t \cos \varphi , t \sin \varphi )^{\rm{t}}$ in polar coordinates, we have that
\begin{equation}
\label{serie_f}
    \mathbf{f}(y)=\sum_{n \in \mathbb{Z}} \mathbf{f}_n(t)e^{i n \varphi},
\end{equation}
with 
\begin{equation*}
    \mathbf{f}_n(t)=\frac{1}{2\pi}\int_0^{2 \pi} \mathbf{f}(t \cos \varphi , t \sin \varphi )e^{-in \varphi} d \varphi.
\end{equation*}

The unique solution of \eqref{ecuacion} satisfying \eqref{radiacionup} and \eqref{radiacionus} can be written as
\begin{align}
    \mathbf{u}(x)
    &=\int_{\mathbb{R}^2} \Phi(x,y) \mathbf{f}(y) dy
\nonumber
    \\
    &= \int_{|y|<|x|} \Phi(x,y) \mathbf{f}(y) dy
    +\int_{|y|>|x|} \Phi(x,y) \mathbf{f}(y) dy
\nonumber
    \\
    &=\mathbf{u}^<(x)+\mathbf{u}^>(x).
\label{u<>}
\end{align}

We will deal first with $\mathbf{u}^<$. Since $|y|<|x|$, we can use Proposition \ref{propo_addition_2D} to write $\Phi(x,y)$ as in \eqref{addition_expansion_2d}. If we split $\mathbf{u}^<$ according to the three terms in the identity \eqref{addition_expansion_2d}, we get
\begin{equation}
\label{u<}
    \mathbf{u}^<=\mathbf{u}^{<,1}+\mathbf{u}^{<,2}<+\mathbf{u}^{<,3},
\end{equation}
where
\begin{equation*}
    \mathbf{u}^{<,j}(x)=\int_{|y|<|x|} \Phi^j(x,y) \mathbf{f}(y) dy,
    \qquad
    j=1,2,3.
\end{equation*}

From here, using polar coordinates, \eqref{Phi1}-\eqref{Phi3} and \eqref{serie_f} and orthogonality, we get
\begin{align}
    \label{uj<}
    \mathbf{u}^{<,j}(x)=\sum_{n \in \mathbb{Z}}\mathbf{u}_n^{<,j}(r)e^{in\theta}, 
    \qquad
    j=1,2,3,
\end{align}
with
\begin{equation}
    \label{un1<}
    \mathbf{u}_n^{<,1}(r)
    =\frac{i}{8}
    \int_0^r t\ H_{n,n}^+(r,t)\ \mathbf{f}_n(t)dt,
\end{equation}
\begin{equation}
    \label{un2<}
    \mathbf{u}_n^{<,2}(r)
    =\frac{1}{16}
    \int_0^r t\ H_{n,n+2}^-(r,t)\ \mathrm{A}\mathbf{f}_{n+2}(t)dt,
\end{equation}
\begin{equation}
    \label{un3<}
    \mathbf{u}_n^{<,3}(r)
    =\frac{1}{16}
    \int_0^r t\ H_{n,n-2}^-(r,t)\ \mathrm{B}\mathbf{f}_{n-2}(t)dt.
\end{equation}

To deal with $\mathbf{u}^>$, we take into account that $\Phi(x,y)=\Phi(y,x)$ (see Remark \ref{r<t}), and thus 
\begin{equation*}
    \mathbf{u}^>(x)
    =\int_{|y|>|x|} \Phi(y,x) \mathbf{f}(y) dy.
\end{equation*}
Now, since we have $\Phi(y,x)$ with $|y|>|x|$, we can use again Proposition \ref{propo_addition_2D} to write $\Phi(y,x)$ as in \eqref{addition_expansion_2d}. And thus, arguing as before, we obtain
\begin{equation}
\label{u>}
    \mathbf{u}^>=\mathbf{u}^{>,1}+\mathbf{u}^{>,2}<+\mathbf{u}^{>,3},
\end{equation}
where
\begin{align}
    \label{uj>}
    \mathbf{u}^{>,j}(x)=\sum_{n \in \mathbb{Z}}\mathbf{u}_n^{>,j}(r)e^{in\theta}, 
    \qquad
    j=1,2,3,
\end{align}
with
\begin{equation}
    \label{un1>}
    \mathbf{u}_n^{>,1}(r)
    =\frac{i}{8}
    \int_r^\infty t\ H_{n,n}^+(t,r)\ \mathbf{f}_n(t)dt,
\end{equation}
\begin{equation}
    \label{un2>}
    \mathbf{u}_n^{>,2}(r)
    =\frac{1}{16}
    \int_r^\infty t\ H_{n+2,n}^-(t,r)\ \mathrm{A}\mathbf{f}_{n+2}(t)dt,
\end{equation}
\begin{equation}
    \label{un3>}
    \mathbf{u}_n^{>,3}(r)
    =\frac{1}{16}
    \int_r^\infty t\ H_{n-2,n}^-(t,r)\ \mathrm{B}\mathbf{f}_{n-2}(t)dt.
\end{equation}
We want to point out that to get \eqref{un1>}-\eqref{un3>} we have made a chnage of variables in the summation index $n$ appearing in \eqref{uj>} and used \eqref{cambio_signo}.

From \eqref{u<>}, \eqref{u<} and \eqref{u>}, to prove \eqref{estimacion_u} it is enough to prove the estimate for $\mathbf{u}^{<,j}$ and $\mathbf{u}^{>,j}$, with $j=1,2,3$. 

We start with $\mathbf{u}^{>,j}$ and show that
\begin{equation}
	\label{estimacion_u>j}
	\int_{\mathbb{R}^2}|\mathbf{u}^{>,j}(x)|^2V(x)dx
	\le c\,|||V|||^2
	\int_{\mathbb{R}^2}|\mathbf{f}(x)|^2V^{-1}(x)dx, \qquad j=1,2,3.
\end{equation}

Writing \eqref{estimacion_u>j} in polar coordinates, using \eqref{uj>} and \eqref{serie_f}, the fact that $V$ is a radial function, and using orthogonality, we see that it is enough to prove the following estimates:
\begin{equation}
	\label{estimacion_u_n>j}
   \sum_{n \in \mathbb{Z}} \int_0^\infty  |\mathbf{u}_n^{>,j}(r)|^2 rV(r) dr
    \le c\,|||V|||^2
    \sum_{n \in \mathbb{Z}}\int_0^\infty  |\mathbf{f}_n(r)|^2 rV^{-1}(r) dr,
    \qquad j=1,2,3,
\end{equation}
where $\mathbf{u}_n^{>,j}$ are given in \eqref{un1>}-\eqref{un3>}.

Arguing in a similar way, the proof of \eqref{estimacion_u>j} for $\mathbf{u}^{<,j}$ can be reduced to proving \eqref{estimacion_u_n>j} but for $\mathbf{u}_n^{<,j}$, which are given in \eqref{un1<}-\eqref{un3<}.

To prove \eqref{estimacion_u_n>j} for $j=1$, we use \eqref{un1>} and the Cauchy-Schwarz inequality, to write
\begin{align*}
    \int_0^\infty  |\mathbf{u}_n^{>,1}(r)|^2 rV(r) dr
    &=\frac{1}{8}\,\int_0^\infty  \left| \int_r^\infty t\ H_{n,n}^+(t,r)\ \mathbf{f}_n(t)dt\right|^2 rV(r) dr
    \\
    &\le
    \int_0^\infty  |\mathbf{f}_n(t)|^2 tV^{-1}(t) dt
    \int_0^\infty \int_r^\infty |\ H_{n,n}^+(t,r)|^2 tV(t) dt\,rV(r) dr,
\end{align*}
reducing matters to establishing
\begin{equation*}
    \int_0^\infty \int_r^\infty |\ H_{n,n}^+(t,r)|^2 tV(t) dt\,rV(r) dr
    \le c\,|||V|||^2,
\end{equation*}
for a constant $c$ independent of $n$.

By \eqref{H+}, this is equivalent to proving
\begin{equation*}
    \int_0^\infty \int_r^\infty |H_n^{(1)}(kt)|^2|J_n(kr)|^2 tV(t) dt\,rV(r) dr
    \le c\,|||V|||^2,
    \qquad k=k_p,k_s.
\end{equation*}

By a change of variables, one can see that this inequality is equivalent to the following:
\begin{equation*}
    \int_0^\infty \int_r^\infty |H_n^{(1)}(t)|^2|J_n(r)|^2 tV_k(t) dt\,rV_k(r) dr
    \le c\,k^2\,|||V|||^2,
\end{equation*}
where $V_k(r)=V(r/k)$. And from \eqref{normaV_a}, this is equivalent to showing
\begin{equation*}
    \int_0^\infty \int_r^\infty |H_n^{(1)}(t)|^2|J_n(r)|^2 tV(t) dt\,rV(r) dr
    \le c\,|||V|||^2.
\end{equation*}
But this estimate was proved in \cite{BRV1997}; more precisely, it is the estimate (2.23) on page 368 of \cite{BRV1997}.

The estimate \eqref{estimacion_u_n>j} for $\mathbf{u}_n^{<,1}$ can be proved analogously, using \eqref{un1<} instead of \eqref{un1>}. In this case, it is necessary to change the order of integration.

To prove \eqref{estimacion_u_n>j} for $j=2$ and $j=3$, we can argue in a similar way, using \eqref{un2>} and \eqref{un3>} instead of \eqref{un1>}, respectively, and we see that it is enough to prove 
\begin{equation}
	\label{estimacion_H-}
    \int_0^\infty \int_r^\infty |\ H_{n,m}^-(t,r)|^2 tV(t) dt\,rV(r) dr
    \le c\,|||V|||^2,\qquad m=n+2, n-2.
\end{equation}
Here we have used that 
$|\mathrm{A}\mathbf{f}_{n+2}|\lesssim |\mathbf{f}_{n+2}|$ and 
$|\mathrm{B}\mathbf{f}_{n-2}|\lesssim |\mathbf{f}_{n-2}|$.

By \eqref{H-}, this is equivalent to proving
\begin{equation*}
    \int_0^\infty \int_r^\infty \left|k_p^2H_m^{(1)}(k_pt)J_n(k_pr)-k_s^2H_m^{(1)}(k_st)J_n(k_sr)\right|^2 tV(t) dt\,rV(r) dr 
    \le c\,|||V|||^2.
\end{equation*}

For the case $k_s>k_p$, we make the change of variables $k_pr=r'$ and $k_pt=t'$ to get (writing $t$ and $r$ instead of $t'$ and $r'$ respectively)
\begin{equation*}
    \int_0^\infty \int_r^\infty \left|H_m^{(1)}(t)J_n(r)-aH_m^{(1)}(at)J_n(ar)\right|^2 tV_{k_p}(t)\,rV_{k_p}(r) dr  dt
    \le c\,k_p^2\, |||V|||^2,
\end{equation*}
where $a=k_s/k_p>1$ and $V_{k_p}(r)=V(r/k_p)$. 
As before, from \eqref{normaV_a}, this inequality is equivalent to 
\begin{equation}
\label{est_dificil}
    \mathrm{I}:=\int_0^\infty \int_r^\infty \left|H_m^{(1)}(t)J_n(r)-aH_m^{(1)}(at)J_n(ar)\right|^2 tV(t) dt\,rV(r) dr
    \le c\, |||V|||^2,
\end{equation}
with $a>1$.

For the case $k_s<k_p$, arguing in the same way but exchanging the roles of $k_p$ and $k_s$, we arrive at the same inequality.

Arguing in a similar way, the proof of \eqref{estimacion_u_n>j} for $\mathbf{u}_n^{<,2}$ or $\mathbf{u}_n^{<,3}$ can be reduced to prove the inequality \eqref{est_dificil}. In this case, one has to use \eqref{un2<} or \eqref{un3<} instead of \eqref{un1>}, and it is necessary to change the order of integration.

So it suffices to prove \eqref{est_dificil} for $m=n+2$ and $m=n-2$ with $n\in\mathbb{Z}$. 
Actually, from \eqref{cambio_signo} it is enough to do this with $n$ a non-negative integer.
Furthermore, from Lemmas \ref{lema_BRV} and \ref{mu_peque} it is sufficient to do so with $n$ sufficiently large.

The proof of \eqref{est_dificil} for $m=n-2$ is much easier than for $m=n+2$ because it is not necessary to take advantage of any cancellation in the integrand. 

Using the triangle inequality, making a change of variables, and employing \eqref{normaV_a}, the proof of \eqref{est_dificil} for $m=n-2$ is reduced to proving that
\begin{equation*}
        \int_0^\infty\int_r^\infty
        |H^{(1)}_\mu(t)|^2
        V(t)t\,dt\,|J_{\mu+2}(r)|^2V(r)r\,dr
        \lesssim |||V|||^2,
\end{equation*}
for $\mu$ large enough.
And this is deduced from \eqref{est_buena} and \eqref{est_BRV}.

\begin{figure}
    \centering
    \begin{tikzpicture}[scale=0.8]
  \draw[->] (0,0) -- (7,0) node[right] {$t$};
  \draw[->] (0,0) -- (0,6.5) node[above] {$r$};

  \draw (0,0) -- (6,6);
  \draw (1,0) -- (1,1);
  \draw (3,2.75) -- (3,3);
  \draw (1,0.75) -- (3,2.75);
  \draw (3,2.75) -- (5,2.75);
  \draw (5,0) -- (5,5);
  
  \fill (0,0) circle (2pt) node[below left] {\footnotesize$A$};
  \fill (1,1) circle (2pt) node[above left] {\footnotesize$B$};
  \fill (3,3) circle (2pt) node[above left] {\footnotesize$C$};
  \fill (5,5) circle (2pt) node[above left] {\footnotesize$D$};
  \fill (5,2.75) circle (2pt) node[below right] {\footnotesize$E$};
  \fill (3,2.75) circle (2pt) node[below right] {\footnotesize$F$};
  \fill (1,0.75) circle (2pt) node[below right] {\footnotesize$G$};
  \fill (1,0) circle (2pt) node[below right] {\footnotesize$H$};
  \fill (5,0) circle (2pt) node[below right] {\footnotesize$I$};

  \node at (6, 2.5) {\footnotesize$D_5$};
  \node at (4, 3.25) {\footnotesize$D_4$};
  \node at (3, 1) {\footnotesize$D_3$};
  \node at (1.9, 2.5) {\footnotesize$D_2$};
  \draw[->] (1.9, 2.3) -- (1.9, 1.75);
  \node at (0.75, 0.25) {\footnotesize$D_1$};
  
  \node at (11,3.5) {\footnotesize{$
  	\begin{array}{l}
  	A=(0,0)
  	\\[0.25ex]
   \displaystyle
  	B=\left(1+\frac{1}{2a},1+\frac{1}{2a}\right)
   	\\[1.75ex]
  	\displaystyle
  	C=\left(\frac{\mu}{2a},\frac{\mu}{2a}\right)
  	\\[1.5ex]
  	D=\left(\mu+2-(\mu+2)^{\frac{1}{3}},\mu+2-(\mu+2)^{\frac{1}{3}}\right)
  	\\[1ex]
  	\displaystyle
  	E=\left(\mu+2-(\mu+2)^{\frac{1}{3}},\frac{\mu}{2a}-\frac{1}{a}\right)
  	\\[2ex]
  	\displaystyle
  	F=\left(\frac{\mu}{2a},\frac{\mu}{2a}-\frac{1}{a}\right)
  	\\[2ex]
  	\displaystyle
  	G=\left(1+\frac{1}{2a},1-\frac{1}{2a}\right)
  	\\[2ex]
  	\displaystyle
  	H=\left(1+\frac{1}{2a},0\right)
  	\\[1ex]
  	\displaystyle
  	I=\left(\mu+2-(\mu+2)^{\frac{1}{3}},0\right)
  	\end{array}
  	$}};
\end{tikzpicture}
    \caption{The integration domain $\{(t,r):r>0,t>r\}$ has been divided into five regions $D_1,\ldots,D_5$.}
    \label{dibu}
\end{figure}
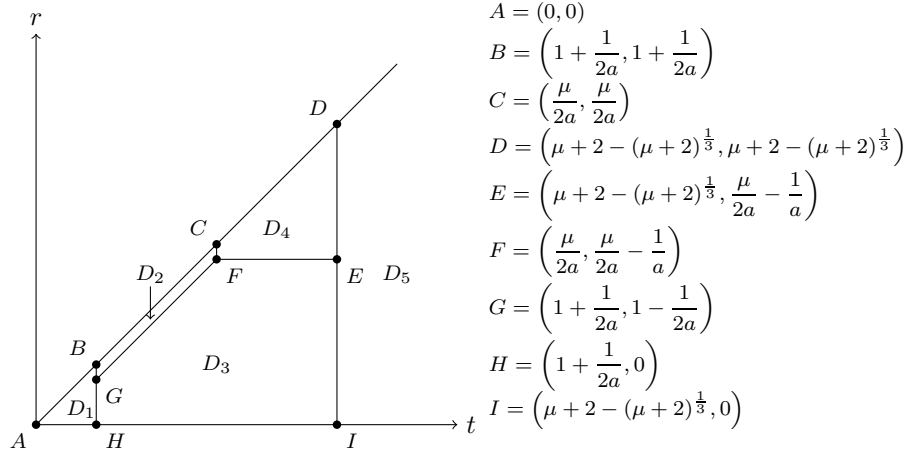

In order to prove \eqref{est_dificil} for $m=n+2$, we split $\mathrm{I}$ into five integrals,
\begin{equation*}
    \mathrm{I}\le \sum_{m=1}^{5}\mathrm{I}_m,
\end{equation*}
where
\begin{equation*}
    \mathrm{I}_m=\iint_{D_m}
    \left|H^{(1)}_{\mu+2}(t)J_\mu(r)-a^2H^{(1)}_{\mu+2}(at)J_\mu(ar)\right|^2V(t)tV(r)r\,drdt,
\end{equation*}
with (see Figure \ref{dibu})
    \begin{equation*}
     D_1=\{(r,t):0< t<1+\tfrac{1}{2a}, 0<r< t\},   
    \end{equation*}
    \begin{equation*}
     D_2=\{(r,t):1+\tfrac{1}{2a}< t<\tfrac{\mu}{2a}, t-\tfrac{1}{a}< r<t\},   
    \end{equation*}
$D_3=D_3^a\cup D_3^b$ with
    \begin{equation*}
     D_3^a=\{(r,t):1+\tfrac{1}{2a}< t<\tfrac{\mu}{2a}, 0<r<t-\tfrac{1}{a}\}
    \end{equation*}
and
    \begin{equation*}
     D_3^b=\{(r,t):\tfrac{\mu}{2a}< t<\mu+2-(\mu+2)^{\frac{1}{3}}, 0<r<\tfrac{\mu}{2a}-\tfrac{1}{a}\},   
    \end{equation*}
    \begin{equation*}
     D_4=\{(r,t):\tfrac{\mu}{2a}<t<\mu+2-(\mu+2)^{\frac{1}{3}}, \tfrac{\mu}{2a}-\tfrac{1}{a}<r<t\} 
    \end{equation*}
and
    \begin{equation*}
     D_5=\{(r,t):t>\mu+2-(\mu+2)^{\frac{1}{3}},0< r<t\}.  
    \end{equation*}

The estimte for $\mathrm{I}_1$ and $\mathrm{I}_2$ follows form Lemma \ref{origen} and \ref{banda}, respectively. 
As we will see later,
the key point in obtaining the estimate in these cases is to take advantage of the cancellation of the integrand (see Remark \ref{cancelacion} below). This is not necessary in the other cases, so  we use the triangle inequality to write
\begin{equation*}
    \mathrm{I}_m\lesssim \mathrm{I}_m^1+\mathrm{I}_m^2,
    \qquad m=3,4,5,
\end{equation*}
where
\begin{equation*}
\begin{split}
    \mathrm{I}_m^1
    &:=
    \iint_{D_m}
    |J_\mu(r)|^2|H^{(1)}_{\mu+2}(t)|^2V(t)tV(r)r\,drdt,\\
    \mathrm{I}_m^2
    &:=
    a^4\iint_{D_m}
    |J_\mu(ar)|^2|H^{(1)}_{\mu+2}(at)|^2V(t)tV(r)r\,drdt.
\end{split}
\end{equation*}
Making a change of variables, we see that
\begin{equation*}
    \mathrm{I}_m^2
    =
    \iint_{D'_m}
    |J_\mu(r')|^2|H^{(1)}_{\mu+2}(t')|^2V_a(t')t'V_a(r')r'\,dr'dt',    
\end{equation*}
where $V_a$ is defined by \eqref{V_a} and 
\begin{equation*}
    D'_m=\left\{(r',t'):\left(\frac{r'}{a},\frac{t'}{a}\right)\in D_m\right\}.
\end{equation*}
From here, one can check that the results for $\mathrm{I}_m^2$ follow from the results for $\mathrm{I}_m^1$ and \eqref{normaV_a}.

Finally, the estimates for $\mathrm{I}_3^1$, $\mathrm{I}_4^1$ and $\mathrm{I}_5^1$  follow from Lemma \ref{escalones}, \ref{Juan} and \ref{lema_BRV}, respectively (see Figure \ref{dibu_escalones}). This concludes the proof for the 2-d case.

\underline{3-D CASE}. 
Following the proof for the 2-D case, we
consider $\mathbf{f}\in \mathcal{C}_0^{\infty}(\mathbb{R}^3)$. Writing $y=(t, \theta_y, \varphi_y )^t$ in spherical coordinates, we have that
\begin{equation}
\label{serie_f3D}
    \mathbf{f}(y)=\sum_{n =0}^\infty\sum_{m=-n}^{+n} \mathbf{f}^m_n(t)Y_n^m(\theta_y, \varphi_y),
\end{equation}
with
\begin{equation*}
 \mathbf{f}^m_n(t)=\int_{S^2} \mathbf {f} (t , \bar y )\overline{Y^m_n}(\bar y) d\sigma(\bar y),
\end{equation*}
where $\bar{y}=y/|y|$.

As we did in \eqref{u<>}, we split the unique solution of \eqref{ecuacion} satisfying \eqref{radiacionup} and \eqref{radiacionus} into two:
\begin{equation*}
\mathbf{u}(x)=\mathbf{u}^<(x)+\mathbf{u}^>(x),
\end{equation*}
where
\begin{equation*}
   \mathbf{u}^<(x)
   =\int_{|y|<|x|} \Phi(x,y) \mathbf{f}(y) dy,
   \quad
   \quad
   \mathbf{u}^>(x)
   =\int_{|y|>|x|} \Phi(x,y) \mathbf{f}(y) dy,
\end{equation*}
and $\Phi$ is the fundamental matrix solution given in \eqref{fundamental3D}.

Arguing as in \eqref{u<}, if $|y|<|x|$ we can use Theorem \ref{adition_theorem} to write $\Phi(x,y)$ as in \eqref{additionLame}. So we can split $\mathbf{u}^<$ according to the four matrix terms in \eqref{additionLame}:
\begin{equation}
\label{u<3D}
    \mathbf{u}^<=\mathbf{u}^{<,1}+\mathbf{u}^{<,2}+\mathbf{u}^{<,3} + \mathbf{u}^{<,4},
\end{equation}
where
\begin{equation*}
    \mathbf{u}^{<,j}(x)=\int_{|y|<|x|} \Phi^j(x,y) \mathbf{f}(y) dy,
    \qquad
    j=1,2,3,4,
\end{equation*}
with $\Phi^1=\Psi,\Phi^2=\Phi_-,\Phi^3=\Phi_0$ and $\Phi^4=\Phi_+.$

From here, using spherical coordinates, \eqref{Phi1_3D}-\eqref{Phi4_3D}, \eqref{serie_f3D}, and orthogonality, we obtain
\begin{equation}
    \label{serie_u_3d}
    \mathbf{u}^{<,j}(x)=\sum_{n =0}^\infty\sum_{m=-n}^{+n}(\mathbf{u}^{<,j})^m_n(r)Y^m_n(\bar x), 
    \qquad
    j=1,2,3,4,
\end{equation}
where $(\mathbf{u}^{<,j})^m_n(r)$ are coefficients written in terms of the corresponding matrices given in \eqref{SP}-\eqref{S+}. These coefficients are uniformly bounded in the corresponding range of $n$ and $m$.  
We can treat the four  matrix terms in  a similar way. To avoid cumbersome repetitions we concentrate on the fourth term given in \eqref{Phi4_3D}. Using orthogonality and making some changes in the summation indicies, we see that
\begin{equation}
\label{u<4_3d}
    (\mathbf{u}^{<,4})_n^m(r)
    =
    \int_0^r h_{n ,n-2}^-(r,t)
    \mathbf{g}^{m}_{n-2}(t)t^2dt,
\end{equation}
with 
\begin{equation}
    \mathbf{g}^{m}_{n-2}=
    \mathrm{A}_{+,n-2}^{m+2}\mathbf{f}^{m+2}_{n-2}+
    \mathrm{B}_{+,n-2}^{m+1}\mathbf{f}^{m+1}_{n-2}+
    \mathrm{C}_{+,n-2}^{m}\mathbf{f}^{m}_{n-2}+
    \mathrm{D}_{+,n-2}^{m-1}\mathbf{f}^{m-1}_{n-2}+
    \mathrm{E}_{+,n-2}^{m-2}\mathbf{f}^{m-2}_{n-2}.
    \label{combinacion_fnm}
\end{equation} 
We would like to point out that the matrix coefficients appearing in the previous identity vanish when $n$ and $m$ are outside the allowed ranges. 

The scalar functions $h_{n ,n-2}^-$ are given in \eqref{h-} in terms of the spherical Bessel functions. 
For convenience, we write them in terms of the standard Bessel and Hankel functions. This gives us (up to a constant)
\begin{align*}
h_{n ,n-2}^-(r,t)
\approx
&(rt)^{-1/2}\left(
k_p^2H^{(1)}_{n+1/2}(k_pr)J_{n-3/2}(k_pt)
-k_s^2H^{(1)}_{n+1/2}(k_sr)J_{n-3/2}(k_s t)
\right)\\
=
&(rt)^{-1/2}H^-_{n+1/2,n-3/2}(r,t),
\end{align*}
where $H^-_{n+1/2,n-3/2}$ is given in \eqref{H-}.

Inserting this in \eqref{u<4_3d} we get
\begin{equation}
\label{u<4nm}
    (\mathbf{u}^{<,4})_n^m(r)
    \approx
    \int_0^r (rt)^{-1/2}H_{n+1/2,n-3/2}^-(r,t)
    \mathbf{g}^{m}_{n-2}(t)t^2dt.
\end{equation}

Following the proof for the 2-d case, matters are reduced to the following estimate:
\begin{equation}
	\label{estimacion_u>j_3d}
	\int_{\mathbb{R}^3}|\mathbf{u}^{<,4}(x)|^2V(x)dx
	\lesssim |||V|||^2
	\int_{\mathbb{R}^3}|\mathbf{f}(x)|^2V^{-1}(x)dx.
\end{equation}
Making the change to spherical coordinates, using \eqref{serie_u_3d} and \eqref{serie_f3D},
and Parseval's identity, we are reduced to prove
\begin{equation*}
    \int_0^\infty r^2 |(\mathbf{u}^{<,4})_n^m(r)|^2V(r)dr
    \le c\, |||V|||^2
    \int_0^\infty r^2 |\mathbf{f}_n^m(r)|^2V^{-1}(r)dr,
\end{equation*}
with $c$ a constant independent of $n$ and $m$.

Using \eqref{u<4nm}, the triangle inequality and Cauchy-Schwarz, since the matrix coefficients appearing in \eqref{combinacion_fnm} are uniformly bounded, we get
\begin{multline*}
\int_0^\infty r^2 |(\mathbf{u}^{<,4})_n^m(r)|^2V(r)dr
\lesssim
\\
\int_0^\infty
\int_0^r
|H^-_{n+\frac{1}{2},n-\frac{3}{2}}(r,t)|^2
tV(t)\,dt\, rV(r)\,dr
\sum_{\ell=-2}^{2}\int_0^rt^2|\mathbf{f}^{m+\ell}_{n-2}(t)|^2V^{-1}(t)dt.
\end{multline*}
From here, changing the order of integration, it is enough to prove 
\begin{equation}
    \int_0^\infty \int_r^\infty |\ H_{\mu,\mu-2}^-(t,r)|^2 tV(t) dt\,rV(r) dr
    \le c\,|||V|||^2,\qquad \mu=n+\frac{1}{2}.
\end{equation}
The proof of this inequality is similar to \eqref{estimacion_H-}.
\hfill $\square$

\begin{remark}
\label{cancelacion}
In order to control the iterated integral $\mathrm{I}$ appearing in \eqref{est_dificil} near the origin, it is important to take advantage of the cancellation of the integrand. This can be seen very clearly taking into account that
$$
H_{n+2}^{(1)}(r) \approx \frac{2^{n+2}}{(n+1)!\,r^{n+2}} ,\;\; r \rightarrow 0, 
\qquad
J_n(t) \approx\frac{t^n}{2^n(n-1)!} ,\;\; t \rightarrow 0,
$$
and then
 $$H_{n+2}(r)J_n(t) \approx \left(  \frac{t}{r} \right)^n
\frac{1}{r^2},$$ 
and this produces a singularity at the origin which is not integrable. 
\end{remark}
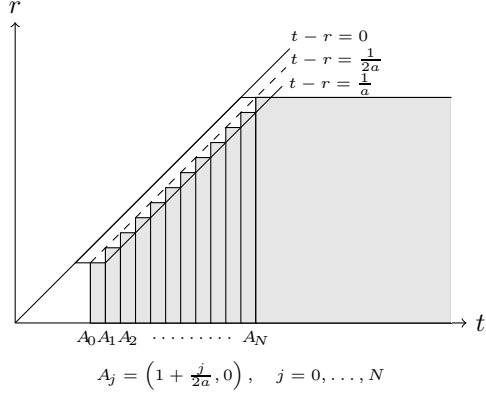
\begin{figure}
    \centering
    \begin{tikzpicture}
  \draw[->] (0,0) -- (6,0) node[right] {$t$};
  \draw[->] (0,0) -- (0,4) node[above] {$r$};
  
    \coordinate (A) at (0.8,0.8);
    \coordinate (B) at (3,3);
    \coordinate (C) at (3.4,3);
    \coordinate (D) at (1.2,0.8);

    \draw (A) -- (B) -- (C) -- (D) -- cycle;


    \coordinate (A1) at (1,0.8);
    \coordinate (B1) at (1.2,0.8);
    \coordinate (C1) at (1.2,0);
    \coordinate (D1) at (1,0);

    \draw (A1) -- (B1) -- (C1) -- (D1) -- cycle;

    \filldraw[fill=gray!20, draw=black] (A1) -- (B1) -- (C1) -- (D1) -- cycle;

    \coordinate (A2) at (1.2,1);
    \coordinate (B2) at (1.4,1);
    \coordinate (C2) at (1.4,0);
    \coordinate (D2) at (1.2,0);

    \draw (A2) -- (B2) -- (C2) -- (D2) -- cycle;

    \filldraw[fill=gray!20, draw=black] (A2) -- (B2) -- (C2) -- (D2) -- cycle;

    \coordinate (A3) at (1.4,1.2);
    \coordinate (B3) at (1.6,1.2);
    \coordinate (C3) at (1.6,0);
    \coordinate (D3) at (1.4,0);

    \draw (A3) -- (B3) -- (C3) -- (D3) -- cycle;

    \filldraw[fill=gray!20, draw=black] (A3) -- (B3) -- (C3) -- (D3) -- cycle;

    \coordinate (A4) at (1.6,1.4);
    \coordinate (B4) at (1.8,1.4);
    \coordinate (C4) at (1.8,0);
    \coordinate (D4) at (1.6,0);

    \draw (A4) -- (B4) -- (C4) -- (D4) -- cycle;

    \filldraw[fill=gray!20, draw=black] (A4) -- (B4) -- (C4) -- (D4) -- cycle;

    \coordinate (A5) at (1.8,1.6);
    \coordinate (B5) at (2,1.6);
    \coordinate (C5) at (2,0);
    \coordinate (D5) at (1.8,0);

    \draw (A5) -- (B5) -- (C5) -- (D5) -- cycle;

    \filldraw[fill=gray!20, draw=black] (A5) -- (B5) -- (C5) -- (D5) -- cycle;

    \coordinate (A6) at (2,1.8);
    \coordinate (B6) at (2.2,1.8);
    \coordinate (C6) at (2.2,0);
    \coordinate (D6) at (2,0);

    \draw (A6) -- (B6) -- (C6) -- (D6) -- cycle;

    \filldraw[fill=gray!20, draw=black] (A6) -- (B6) -- (C6) -- (D6) -- cycle;

    \coordinate (A7) at (2.2,2);
    \coordinate (B7) at (2.4,2);
    \coordinate (C7) at (2.4,0);
    \coordinate (D7) at (2.2,0);

    \draw (A7) -- (B7) -- (C7) -- (D7) -- cycle;

    \filldraw[fill=gray!20, draw=black] (A7) -- (B7) -- (C7) -- (D7) -- cycle;

    \coordinate (A8) at (2.4,2.2);
    \coordinate (B8) at (2.6,2.2);
    \coordinate (C8) at (2.6,0);
    \coordinate (D8) at (2.4,0);

    \draw (A8) -- (B8) -- (C8) -- (D8) -- cycle;

    \filldraw[fill=gray!20, draw=black] (A8) -- (B8) -- (C8) -- (D8) -- cycle;

    \coordinate (A9) at (2.6,2.4);
    \coordinate (B9) at (2.8,2.4);
    \coordinate (C9) at (2.8,0);
    \coordinate (D9) at (2.6,0);

    \draw (A9) -- (B9) -- (C9) -- (D9) -- cycle;

    \filldraw[fill=gray!20, draw=black] (A9) -- (B9) -- (C9) -- (D9) -- cycle;

    \coordinate (A10) at (2.8,2.6);
    \coordinate (B10) at (3,2.6);
    \coordinate (C10) at (3,0);
    \coordinate (D10) at (2.8,0);

    \draw (A10) -- (B10) -- (C10) -- (D10) -- cycle;

    \filldraw[fill=gray!20, draw=black] (A10) -- (B10) -- (C10) -- (D10) -- cycle;

    \coordinate (A11) at (3,2.8);
    \coordinate (B11) at (3.2,2.8);
    \coordinate (C11) at (3.2,0);
    \coordinate (D11) at (3,0);

    \draw (A11) -- (B11) -- (C11) -- (D11) -- cycle;

    \filldraw[fill=gray!20, draw=black] (A11) -- (B11) -- (C11) -- (D11) -- cycle;
    
    \coordinate (Au) at (3.2,3);
    \coordinate (Bu) at (5.8,3);
    \coordinate (Cu) at (5.8,0);
    \coordinate (Du) at (3.2,0);

    \draw (Au) -- (Bu) -- (Cu) -- (Du) -- cycle;

    \filldraw[fill=gray!20, draw=black] (Au) -- (Bu) -- (Cu) -- (Du) -- cycle;

  \draw (0,0) -- (3.65,3.65);
  \draw[dashed] (1,0.8) -- (3.6,3.4);
  \draw (1.2,0.8) -- (3.55,3.15);
  \draw[white] (5.8,0) -- (5.8,3);

  \node at (0.95,-0.2) {\tiny$A_{\!0}$};
  \node at (1.225,-0.2) {\tiny$A_{\!1}$};
  \node at (1.5,-0.2) {\tiny$A_{\!2}$};
  \node at (2.35,-0.2) {\tiny$\ldots\ldots\ldots$};
  \node at (3.2,-0.2) {\tiny$A_{\!N}$};

  \node at (3,-0.7) {\tiny$A_{\!j}=\left(1+\frac{j}{2a},0\right),\quad j=0,\ldots, N$};

  \node at (4.18,3.8) {\tiny$t-r=0$};
  \node at (4.28,3.5) {\tiny$t-r=\frac{1}{2a}$};
  \node at (4.2,3.18) {\tiny$t-r=\frac{1}{a}$};
\end{tikzpicture}
    \caption{Region covering the integration domains of Lemma \ref{escalones}.}
    \label{dibu_escalones}
\end{figure}
\section{Proof of Theorem \ref{teorema_gradiente}}
We will need the following lemma, which can be proved following the proof of Proposition 2(a) in page 245 of \cite{Stein}.

\begin{lemma}\cite[pp. 245, Proposition 2]{Stein}
\label{malo_lema}
    Let $m$ be a bounded function in $C^\infty(\mathbb{R}^d\setminus\{0\})$ that satisfies
    \begin{equation}
        \label{cota_derivadas_m}
        |\partial^\alpha m(\xi)|\le 
        c_\alpha|\xi|^{-(1+|\alpha|)},
    \end{equation}
    for all $\alpha=(\alpha_1,\ldots,\alpha_d)\in\mathbb{N}^d$ with $|\alpha|\le d$.
    And denote by $K$ the distribution whose Fourier transform is $m$. Then, $K$ satisfies
    \begin{equation}
        \label{cota_K}
        |K(x)|\lesssim |x|^{-(d-1)}.
    \end{equation}
\end{lemma}
We will also need the following estimate for a Riesz potential.
\begin{proposition}
    \cite[pp. 365, (2.15)]{BRV1997} 
    \label{fraccionaria_th}
    If $V\in \mathcal{T}(\mathbb{R}^d)$,  then there exists a constant $c>0$ independent of $f$ such that the following estimate holds:
    \begin{equation}
        \label{fraccionaria_est}
        \int_{\mathbb{R}^d}| I_1f(x)|^2V(x)dx\le c\, |||V|||^2\int_{\mathbb{R}^d}|f(x)|^2V^{-1}(x)dx,
    \end{equation}
where $I_1$ is the Riesz potential defined by 
\begin{equation*}
    I_1f(x)=c_d\, |x|^{-(d-1)}\ast f(x).
\end{equation*}
\end{proposition}
We would like to remark
that by duality, estimate \eqref{fraccionaria_est} follows from the following estimate:
\begin{equation*}
         \|V^{1/2}I_{1/2}f\|_2
         \le c\,
         M_{1,2,\infty}(V^{1/2})\|f\|_2,
    \end{equation*}
where $M_{1,2,\infty}(V^{1/2})=\sup_xI_1(V)(x)$, and
which can be proved following the lines of Theorem 7.2 in Chapter 6 in \cite{Sch}.
    
We will also need the following known result for the solution of the Helmholtz equation.
\begin{theorem}\cite[pp. 362, (1.16)]{BRV1997}
    \label{solucion_th}
    Let $f\in C_0^\infty (\mathbb{R}^d)$  and consider $u$, the unique solution to problem \eqref{Helmholtz_eq} satisfying \eqref{SRC}. If $V_1,V_2\in \mathcal{T}_r(\mathbb{R}^d)$, then there exists a constant $c>0$ independent of $f$ and $k$ such that the following a priori estimate holds:
    \begin{equation}
        \label{solucion_est}
        \int_{\mathbb{R}^d}| u(x)|^2V_1(x)dx\le \frac{c}{k^2}\, |||V_1|||\,|||V_2|||\int_{\mathbb{R}^d}|f(x)|^2V_2^{-1}(x)dx.
    \end{equation}
\end{theorem}
\emph{Proof of Theorem \ref{teorema_gradiente}.}
As in the proof of Theorem \ref{principal},
it is sufficient to prove the theorem for the radial majorant of $V$ given in \eqref{mayorante_radial}, which, abusing notation, we will also denote by $V$. As before, we may assume that, without losing generality, $\omega=1$.

Moreover, it suffices to prove 
\begin{equation}
        \label{parcial_est}
        \int_{\mathbb{R}^d}|\partial_{\xi_\ell} \mathbf{u}(x)|^2V(x)dx\le c\, |||V|||^2\int_{\mathbb{R}^d}|\mathbf{f}(x)|^2V^{-1}(x)dx,\qquad \ell=1,\ldots,d,
    \end{equation}
where $\partial_{\xi_\ell} \mathbf{u}$ is defined componentwise.

For convenience, we introduce the Leray projection operator, 
$\mathrm{I}-\mathcal{R}$, where $\mathrm{I}$ is the identity and $\mathcal{R}$ is given by
    \begin{equation*}
        \label{Leray}
        \widehat{\mathcal{R}\mathbf{f}}(\xi)
        =\frac{\xi}{|\xi|}\cdot \widehat{\mathbf{f}}(\xi)\ \frac{\xi}{|\xi|},\qquad\qquad \xi\in\mathbb{R}^d.
    \end{equation*}

Since 
$\nabla div\,\mathbf{u}=\Delta \mathcal{R}u$,
we can see that $\mathcal{R}\mathbf{u}$ and $(\mathrm{I}-\mathcal{R})\mathbf{u}$ are solutions of the following vectorial Helmholtz equations:
\begin{align*}
    \Delta \mathcal{R}\mathbf{u}+\frac{1}{2\mu+\lambda} \mathcal{R}\mathbf{u}&=\frac{1}{2\mu+\lambda} \mathcal{R}\mathbf{f}, 
    \\
    \Delta (\mathrm{I}-\mathcal{R})\mathbf{u}+\frac{1}{\mu} (\mathrm{I}-\mathcal{R})\mathbf{u}&=\frac{1}{\mu} (\mathrm{I}-\mathcal{R})\mathbf{f}.
\end{align*}

From here, since $\mathbf{u}=\mathcal{R}\mathbf{u}+(\mathrm{I}-\mathcal{R})\mathbf{u}$, we can write
\begin{equation*}
    \widehat{\mathbf{u}}(\xi)=\mathrm{M}(\xi)\,\widehat{\mathbf{f}}(\xi),
\end{equation*}
where $\mathrm{M}(\xi)$ is the matrix given by
\begin{equation*}
    \mathrm{M}(\xi)=\mathrm{M}_1(\xi)+\mathrm{M}_2(\xi),
\end{equation*}
with
\begin{align}
    \mathrm{M}_1(\xi)&=
     \frac{\mathrm{I}}{\mu\left(\frac{1}{\mu}-|\xi|^2+i0\right)}, \ \ \ {\rm and}
     \nonumber
     \\
    \mathrm{M}_2(\xi)&=
    \frac{\xi_i\xi_j}{|\xi|^2}
    \left(
    \frac{1}{(2\mu+\lambda)\left(\frac{1}{2\mu+\lambda}-|\xi|^2+i0\right)}
    -
    \frac{1}{\mu\left(\frac{1}{\mu}-|\xi|^2+i0\right)}
    \right).
    \label{M2}
\end{align}

Now we split $\mathbf{u}$ into two as follows:
\begin{equation*}
    \mathbf{u}=\mathbf{u}_1+\mathbf{u}_2,
\end{equation*}
where 
\begin{equation*}
    \widehat{\mathbf{u}}_m(\xi)=\mathrm{M}_m(\xi)\,\widehat{\mathbf{f}}(\xi),
    \qquad m=1,2.
\end{equation*}

The result for $\mathbf{u}_1$ follows from estimate \eqref{l-resolvent}.

In order to prove the result for $\mathbf{u}_2$ we introduce a function $\phi\in C_0^\infty(\mathbb{R}^d)$. For the case $\sqrt{\mu}\le\sqrt{2\mu+\lambda}$, we consider $\phi$ supported in 
$B\left(0,\frac{2}{\sqrt{\mu}}\right)\setminus B\left(0,\frac{1}{2\sqrt{2\mu+\lambda}}\right)$ 
such that $\phi(\xi)=1$ if
$\xi\in B\left(0,\frac{3}{2\sqrt{\mu}}\right)\setminus B\left(0,\frac{3}{4\sqrt{2\mu+\lambda}}\right)$.
For the case $\sqrt{2\mu+\lambda}\le \sqrt{\mu}$, we simply interchange the roles between $\mu$ and $2\mu+\lambda$ in the conditions that must be satisfied by $\phi$.

Now we split the multiplier $\mathrm{M}_2$ into two as follows:
\begin{equation*}
    \mathrm{M}_2=(1-\phi)\,\mathrm{M}_2+\phi\, \mathrm{M}_2,    
\end{equation*}
and we write
\begin{equation*}
    (\partial_{\xi_\ell} \mathbf{u}_2)\widehat{\phantom{u}}(\xi)
    =
    i\xi_\ell\, 
    \mathrm{M}_2(\xi)\,\widehat{\mathbf{f}}(\xi)
    =
    \mathrm{M}_2^1(\xi)\,\widehat{\mathbf{f}}(\xi)
    +
    \mathrm{M}_2^2(\xi)\,\widehat{\mathbf{f}}(\xi),
\end{equation*}
where
\begin{equation}
\label{M12M22}
    \mathrm{M}_2^1(\xi)=i\xi_\ell\, (1-\phi(\xi))\,\mathrm{M}_2(\xi),
    \qquad
    \text{ and }
    \qquad
    \mathrm{M}_2^2(\xi)=i\xi_\ell\, \phi(\xi)\,\mathrm{M}_2(\xi).
\end{equation}

From \eqref{M2}, we see that we can write
\begin{equation}
\label{M2_b}
    \mathrm{M}_2(\xi)=
    \frac{(\mu+\lambda)\,\xi_i\xi_j}{\mu(2\mu+\lambda)\left(\frac{1}{\mu}-|\xi|^2+i0\right)\left(\frac{1}{2\mu+\lambda}-|\xi|^2+i0\right)}.
\end{equation}
From this, we see that
each component of the matrix multiplier $\mathrm{M}_2^1$ satisfies \eqref{cota_derivadas_m}, so we can apply Lemma \ref{malo_lema} and conclude that each component of the matrix distribution $\mathrm{K}_1$ whose Fourier transform is $\mathrm{M}_2^1$ satisfies \eqref{cota_K}.

On the other hand, from \eqref{M2} we write
\begin{equation*}
    \mathrm{M}_2^2(\xi)=
    \left(
    \frac{1}{(2\mu+\lambda)\left(\frac{1}{2\mu+\lambda}-|\xi|^2+i0\right)}
    -
    \frac{1}{\mu\left(\frac{1}{\mu}-|\xi|^2+i0\right)}
    \right)\,\mathrm{N}(\xi),
\end{equation*}
where
\begin{equation*}
    \mathrm{N}(\xi)=i\xi_\ell\phi(\xi)\,\frac{\xi_i\xi_j}{|\xi|^2}.
\end{equation*}
Note that each component of the matrix $\mathrm{N}$ is a $C_0^\infty$ function. We denote them by $n_{ij}$. Hence
the functions $k_{ij}$ whose Fourier transforms are $n_{ij}$ are in the Schwartz class. We consider the matrix $\mathrm{K}:=(k_{ij})$.

Taking all of this into account, we can write
\begin{equation*}
    \partial_{\xi_j} \mathbf{u}_2(x)
    =\mathrm{K}_1\ast \mathbf{f}(x)
    +\mathrm{K}\ast
    \left(
    \mathbf{v}-\mathbf{w}
    \right)(x),
\end{equation*}
where $\mathbf{v}$ and $\mathbf{w}$ are solutions of the following vectorial Helmholtz equations:
\begin{align*}
    \Delta \mathbf{v} +\frac{1}{2\mu+\lambda}\mathbf{v}&=\frac{1}{2\mu+\lambda}\mathbf{f},
    \\
    \Delta \mathbf{w} +\frac{1}{\mu}\mathbf{w}&=\frac{1}{\mu}\mathbf{f}.
\end{align*}

The result for the first term, $K_1\ast \mathbf{f}$, follows from Proposition \ref{fraccionaria_th} and the fact that $K_1$ satisfies \eqref{cota_K}.

In order to prove the result for the second term, we use the Cauchy-Schwarz inequality to write
\begin{align}
    \int_{\mathbb{R}^d}|K\ast \mathbf{v}(x)|^2V(x)dx
    &=
    \int_{\mathbb{R}^d}\left|\int_{\mathbb{R}^d} K(y)\mathbf{v}(x-y)dy\right|^2V(x)dx
    \nonumber
    \\
    &\le 
    \|K\|_{L^1(\mathbb{R}^d)}
    \int_{\mathbb{R}^d}\int_{\mathbb{R}^d}|K(y)||\mathbf{v}(x-y)|^2dyV(x)dx
    \nonumber
    \\
    &\le c\,
    \int_{\mathbb{R}^d}|\mathbf{v}(y)|^2\int_{\mathbb{R}^d}|K(x-y)|V(x)dxdy
    \nonumber
    \\
    &\le c\,
    \int_{\mathbb{R}^d}|\mathbf{v}(y)|^2W(y)dy,
    \label{final}
\end{align}
where $W=\widetilde{K}\ast V$ with 
\begin{equation*}
   \widetilde{K}(y)=\sup_{x\in S_{|y|}} |K(x)|. 
\end{equation*} 

On the other hand, for any $x\in\mathbb{R}^d$ and $\omega\in S^{d-1}$, changing the order of integration, we see that
\begin{align*}
    \int_{\mathbb{R}}W(x+t\omega)\,dt
    &=
    \int_{\mathbb{R}^d}\widetilde{K}(y)\int_{\mathbb{R}}V(x+t\omega-y)\,dt\,dy
    \\
    &\le
    \|\widetilde{K}\|_{L^1}|||V|||
    \\
    &\lesssim
    |||V|||.
\end{align*}
Therefore, since $K$ is in the Schwartz class, $W\in \mathcal{T}_{rad}(\mathbb{R}^n)$ and
$|||W|||\lesssim |||V|||$. Using Theorem \ref{solucion_th} with $V_1=W$ and $V_2=V$ in \eqref{final}, the result follows.

The same holds for $\mathbf{w}$.
\hfill $\square$
\section{Proof of the key lemmas.}
\label{demos_lemas}
In this section we will need several known results for Bessel and Hankel functions. 
Firstly, we present the ascending series expansions of the Bessel functions of the first and second kind, which can be found in \cite[p. 360]{AS1967}.
\begin{lemma}\cite[Identity 9.1.10]{AS1967}
    \label{Lema_seriesJ}
    Let $\mu\in \mathbb{R}$ and $r\ge 0$, the ascending series expansion of the first Bessel function, denoted by $J_\mu$, is given by
    \begin{equation}
        \label{serie_J}
        J_\mu(r)=\left(\frac{r}{2}\right)^\mu\sum_{k=0}^\infty\frac{\left(-\frac{r^2}{4}\right)^k}{k!\,\Gamma(\mu+k+1)}.
    \end{equation}
\end{lemma}
\begin{lemma}\cite[Identity 9.1.11]{AS1967}
    \label{Lema_seriesY}
    Let $n$ be a non-negative integer and $r> 0$, the ascending series expansion of the second Bessel function, denoted by $Y_n$, is given by
    \begin{equation}
        \label{serie_Y_n}
        Y_n(r)=-\frac{\left(\frac{r}{2}\right)^{-n}}{\pi}\sum_{k=0}^{n-1}\frac{(n-k-1)!}{k!}\,\left(\frac{r^2}{4}\right)^{k}+O_n(r),
    \end{equation}
    where
    \begin{equation}
        \label{On}
        O_n(r)=\frac{2}{\pi}\,\ln \left(\frac{r}{2}\right)J_n(r)-
        \frac{\left(\frac{r}{2}\right)^{n}}{\pi}\sum_{k=0}^\infty(\psi(k+1)+\psi(n+k+1))\,\frac{\left(\frac{r^2}{4}\right)^k}{k!\,(n+k)!},
    \end{equation}
    with $\psi(1)=-\gamma$ where $\gamma$ is the Euler constant, and
    \begin{equation*}
        \psi(n)=-\gamma+\sum_{k=1}^{n-1}\frac{1}{k}.
    \end{equation*}
\end{lemma}
\begin{remark}
    From the definition of the second Bessel function of fractional order given in \cite[p.64]{Watson} we have that
    $Y_{n+\frac{1}{2}}(r)=(-1)^{n+1}J_{-\left(n+\frac{1}{2}\right)}(r)$. Thus, using the expansion \eqref{serie_J} we get the following ascending series expansion:
    \begin{equation}
    \label{serie_Y_n+1/2}
        Y_{n+\frac{1}{2}}(r)=\left(\frac{r}{2}\right)^{-\left(n+\frac{1}{2}\right)}\sum_{k=0}^\infty\frac{\left(-\frac{r^2}{4}\right)^k}{k!\,\Gamma(-\left(n+\frac{1}{2}\right)+k+1)}.    
    \end{equation}
\end{remark}
\begin{remark}
    \label{origen_Lema}
    From \eqref{serie_J}, we have that
    for $\mu\ge 0$, if $0\le r<2\sqrt{\mu}$, then 
    \begin{equation}
    \label{origen_J}
    J_\mu(r)=\left(\frac{r}{2}\right)^\mu\frac{1}{\Gamma(\mu+1)}
    +
    O\left(\left(\frac{r}{2}\right)^{\mu+2}\frac{1}{\Gamma(\mu+2)}\right),
    \qquad
    \mu\rightarrow \infty.
    \end{equation}
    Moreover, from \eqref{serie_Y_n} and \eqref{serie_Y_n+1/2}, for $\mu=n$ or $\mu=n+\frac{1}{2}$ with $n$ a non-negative integer, if $0< r<2\sqrt{\mu}$, then 
    \begin{equation}
    \label{origen_Y}
    Y_\mu(r)=-\frac{i}{\pi}\left(\frac{r}{2}\right)^{-\mu}\Gamma(\mu)
    +
    O\left(\left(\frac{r}{2}\right)^{-\mu+2}\Gamma(\mu-1)\right),
    \qquad
    \mu\rightarrow \infty.    
    \end{equation}
    We would like to point out that to obtain \eqref{origen_Y} in the case $\mu=n+\frac{1}{2}$, we have used that
    \begin{equation*}
        \Gamma\left(n+\frac{1}{2}\right)\Gamma\left(1-\left(n+\frac{1}{2}\right)\right)=(-1)^n\pi.
    \end{equation*}
\end{remark}

\bigskip
\textit{Proof of Lemma \ref{mu_peque}.}
For $a\ge1$ fixed we denote by $\mathrm{I}$ the iterated integral in \eqref{est_mu_peque}.
First of all, if we denote
\begin{equation}
\label{D_nu_mu}
     D_{\nu,\mu,a}(r,t):= \left|H^{(1)}_\nu(t)J_{\mu}(r)-a^2H^{(1)}_\nu(at)J_{\mu}(ar)\right|,
\end{equation}
we have $\mathrm{I}\le \mathrm{I}_1+\mathrm{I}_2$, where
\begin{align*}
    \mathrm{I}_1&=
    \int_0^{\nu-\nu^{\frac{1}{3}}}\int_r^{\nu-\nu^{\frac{1}{3}}}
    |D_{\nu,\mu,a}(r,t)|^2
    V(t)t\,dt\,V(r)r\,dr,
    \\
    \mathrm{I}_2&=
    \int_0^\infty\int_{\nu-\nu^{\frac{1}{3}}}^\infty
    |D_{\nu,\mu,a}(r,t)|^2V(t)t\,dt\,V(r)r\,dr.
\end{align*}
The estimate for $\mathrm{I}_2$ is easier than for $\mathrm{I}_1$.
Using the triangle inequality and making a change of variables, since $a\ge 1$, we have
\begin{equation*}
\begin{split}
    \mathrm{I}_2\lesssim &
   \int_0^\infty\left|J_\mu(r)\right|^2V(r)r\,dr
   \int_{\nu-\nu^{\frac{1}{3}}}^\infty\left|H^{(1)}_\nu(t)\right|^2V(t)t\,dt\\
   &+\int_0^\infty\left|J_\mu(r)\right|^2V_a(r)r\,dr
   \int_{\nu-\nu^{\frac{1}{3}}}^\infty\left|H^{(1)}_\nu(t)\right|^2V_a(t)t\,dt,
\end{split}
\end{equation*}
where $V_a$ is defined by \eqref{V_a}. 
From here, the estimate for $\mathrm{I}_2$ follows from Lemma \ref{lema_BRV} and identity \eqref{normaV_a}.

In order to prove the estimate for $\mathrm{I}_1$, it is necessary to take advantage of the cancellation in \eqref{D_nu_mu}.
But first, since $H_{\nu}^{(1)}=J_{\nu}+iY_{\nu}$, using the triangle inequality, we split $\mathrm{I}_1$ into two
\begin{equation*}
    \mathrm{I}_1\lesssim \mathrm{I}_1^1+\mathrm{I}_1^2,
\end{equation*}
where
\begin{equation*}
    \mathrm{I}_1^m=\int_0^{\nu-\nu^{\frac{1}{3}}}\int_r^{\nu-\nu^{\frac{1}{3}}}
        \left|D_{\nu,\mu,a}^m(r,t)\right|^2
        V(t)t\,dt\,V(r)r\,dr,
        \qquad m=1,2,
\end{equation*}
with
\begin{align}
\label{Dm1}
    D_{\nu,\mu,a}^1(r,t)&= \left|J_\nu(t)J_{\mu}(r)-a^2J_\nu(at)J_{\mu}(ar)\right|,\\
\label{Dm2}
    D_{\nu,\mu,a}^2(r,t)&= \left|Y_\nu(t)J_{\mu}(r)-a^2Y_\nu(at)J_{\mu}(ar)\right|.
\end{align}

The estimate for $\mathrm{I}_1^1$ is obtained in a similar way as for $\mathrm{I}_2$.

To finish, we have to prove the estimate for $\mathrm{I}_1^2$. We consider $\nu=\mu+2$ since the case $\nu=\mu-2$ is easier. We also consider $\mu=n$, since the case $\mu=n+\frac{1}{2}$ is similar. Using \eqref{serie_J} and \eqref{serie_Y_n} one can see that
\begin{equation*}
    a^2J_n(ar)Y_{n+2}(at)=-\frac{i}{\pi}\left(\frac{r}{t}\right)^n (n+1)\frac{4}{t^2}
    +
    R_{n,a}(r,t),
\end{equation*}
where $R_{n,a}$ is the remainder term. Some rough calculations gives us the following bound for $R_{n,a}(r,t)$ for $(r,t)$ within the integration region of $\mathrm{I}_1^2$: 
\begin{equation*}
    |R_{n,a}(r,t)|\lesssim \left(\frac{r}{t}\right)^n \left(e^{\frac{ar^2}{4}}+e^{\frac{at^2}{4}}\right)\le c_{n,a}.
\end{equation*}
From here we see that
\begin{equation*}
    D_{n+2,n,a}^2(r,t)= |R_{n,1}(r,t)-R_{n,a}(r,t)|\le |R_{n,1}(r,t)|+|R_{n,a}(r,t)|\lesssim c_{n,a},
\end{equation*}
and therefore, 
\begin{align*}
    \mathrm{I}_1^2
    &\lesssim
    c_{n,a}
    \left(\int_0^\infty V(t)\,dt\right)\left(\int_0^\infty V(r)\,dr\right)
    \\
    &\le c_{n,a} |||V|||^2.
\end{align*}
This concludes the proof for the case $\mu=n$. The result for $\mu=n+\frac{1}{2}$ can be obtained  in a similar way by using \eqref{serie_Y_n+1/2} instead of \eqref{serie_Y_n}.
\hfill $\square$

\bigskip
\textit{Proof of Lemma \ref{origen}.}
Arguing as in the proof of Lemma \ref{mu_peque}, it is enough to prove the estimate for 
\begin{equation*}
    \mathrm{I}=\int_0^{1+\frac{1}{2a}}\int_r^{1+\frac{1}{2a}}
        \left|D_{\mu,a}(r,t)\right|^2
        V(t)t\,dt\,V(r)r\,dr,
\end{equation*}
where
\begin{equation}
\label{D_origen}
    D_{\mu,a}(r,t)= \left|J_\mu(r)Y_{\mu+2}(t)-a^2J_\mu(ar)Y_{\mu+2}(at)\right|.
\end{equation}

Using \eqref{origen_J} and \eqref{origen_Y} we can write
\begin{equation*}
    a^2J_\mu(ar)Y_{\mu+2}(at)=-\frac{i}{\pi}\left(\frac{r}{t}\right)^\mu (\mu+1)\frac{4}{t^2}
    +
    O\left(\left(\frac{r}{t}\right)^{\mu}\right),
\end{equation*}
for $\mu>\mu_0$ with $\mu_0=\mu_0(a)$ large enough and $(r,t)$ within the integration region. 
This identity allows us to take advantage of cancellation in \eqref{D_origen} and write
\begin{equation*}
    D_{\mu,a}(r,t)\lesssim 1.
\end{equation*}
From here, the proof concludes as in Lemma \ref{mu_peque}.
\hfill $\square$

Now we present a result that can be found in \cite[p. 118 and 121]{Matviyenko}, where the author gives bounds 
for the error terms of the so-called asymptotic Debye expansions given in \cite[Chap. 10]{Olver}. 
More precisely, considering the identities (32) and (33), and the Theorem 3.3 in \cite{Matviyenko} for the particular case $N=0$, we obtain the following lemma.
\begin{lemma}
    \label{Deby_Lema}
    If $\mu$ is large and $0\le r<\mu$, then the first and the second Bessel functions, denoted by $J_\mu$ and $Y_\mu$ respectively, can be written as follows: 
\begin{align}
    \label{Deby_J}
    J_\mu(r)&=c_0\,\sqrt{\frac{1}{2\pi}}\,
    \frac{e^{-\mu\phi_\mu(r)}}{(\mu^2-r^2)^{1/4}}\left(1+e_{1,\mu}(r)\right),
    \\
    \label{Deby_Y}
    Y_\mu(r)&=-\sqrt{\frac{2}{\pi}}\,
    \frac{e^{\mu\phi_\mu(r)}}{(\mu^2-r^2)^{1/4}}\left(1+e_{2,\mu}(r)\right),
\end{align}
with $c_0$ an absolute constant, 
\begin{equation}
\label{phi}
    \phi_\mu(r)=\alpha_\mu(r)-\tanh\alpha_\mu(r),
\end{equation}
\begin{equation}
    \label{alpha}
    e^{\alpha_\mu(r)}=\frac{\mu+\sqrt{\mu^2-r^2}}{r}
    \text{ and }
    \tanh{\alpha_\mu(r)}=\frac{\sqrt{\mu^2-r^2}}{\mu},
\end{equation}
and where
\begin{align}
    \label{error_J}
    \left|e_{1,\mu}(r)\right|
    \le 
    \frac{2}{3}\left(\frac{\mu^{1/3}}{\mu-r}\right)^{3/2}e^{\frac{2}{3}\left(\frac{\mu^{1/3}}{\mu-r}\right)^{2/3}},
    \\
    \label{error_Y}
    \left|e_{2,\mu}(r)\right|
    \le 
    \frac{2}{3}\left(\frac{\mu^{1/3}}{\mu-r}\right)^{3/2}e^{\frac{2}{3}\left(\frac{\mu^{1/3}}{\mu-r}\right)^{3/2}}.
\end{align}
\end{lemma}
Regarding \eqref{Deby_J}, see also \cite[p.38]{Nemes} (inequality (6.4) and the previous identity).
\begin{remark}
    \label{cota_J_Y}
    Note that from \eqref{error_J} and \eqref{error_Y} we have that the error terms satisfy
    \begin{align}
        \label{errores_cota1}
        \left|e_{\ell,\mu}(r)\right|&\lesssim\mu^{-1}, 
        \qquad\text{ if }0\le r\le\mu/2,\ \ell=1,2,
        \\
        \label{errores_cota2}
        \left|e_{\ell,\mu}(r)\right|&\lesssim 1, 
        \qquad\text{ if }0\le r\le\mu+2-(\mu+2)^{1/3},\ \ell=1,2.
    \end{align}
    Moreover, using \eqref{errores_cota1} in \eqref{Deby_J} and \eqref{Deby_Y}, we obtain
    \begin{align}
    \label{cota_J1}
        \left|J_\mu(r)\right|&\lesssim 
        \frac{e^{-\mu\phi_\mu(r)}}{\mu^{1/2}},
        \qquad\text{ if }0\le r\le\mu/2,
        \\
        \nonumber
        \left|Y_\mu(r)\right|&\lesssim 
        \frac{e^{\mu\phi_\mu(r)}}{\mu^{1/2}},
        \qquad\text{ if }0\le r\le\mu/2.
    \end{align}
    And since $H^{(1)}_\mu=J_\mu+iY_\mu$, we also have that
\begin{equation}
    \label{cota_H1}
        \left|H^{(1)}_\mu(r)\right|\lesssim 
        \frac{e^{\mu\phi_\mu(r)}}{\mu^{1/2}},
        \qquad\text{ if }0\le r\le\mu/2.
\end{equation}
\end{remark}
\textit{Proof of Lemma \ref{banda}.}
For fixed $a\ge1$, we denote by $\mathrm{I}$ the iterated integral in \eqref{est_banda}. 
Arguing as we did in the proof of Lemma \ref{origen_Lema} 
we split $\mathrm{I}$ into two
\begin{equation*}
    \mathrm{I}\lesssim \mathrm{I}_1+\mathrm{I}_2,
\end{equation*}
where
\begin{equation*}
    \mathrm{I}_m=\int_{1-\frac{1}{2a}}^{\frac{\mu}{2a}}\left(\int_r^{r+\frac{1}{a}}
        \left|D_{\mu,a}^m(r,t)\right|^2
        V(t)t\,dt\right)V(r)r\,dr,
        \qquad m=1,2,
\end{equation*}
with $D_{\mu,a}^1$ and $D_{\mu,a}^2$ given in \eqref{Dm1} and \eqref{Dm2} respectively.
Hence, it is enough to prove the estimate for $\mathrm{I}_1$ and $\mathrm{I}_2$.

The estimate for $\mathrm{I}_1$ is obtained in the same way as in the proof of Lemma \ref{origen_Lema}.
In order to prove the estimate for $\mathrm{I}_2$, we use \eqref{Deby_J} and \eqref{Deby_Y} to write
\begin{equation*}
\begin{split}
    D_{\mu,a}^2(r,t)
    \sim &
    \left|
    F_\mu(r,t)(1+e_{1,\mu}(r))(1+e_{2,\mu+2}(t))
    \right.
    \\&-
    \left.
    a^2F_\mu(ar,at)(1+e_{1,\mu}(ar))(1+e_{2,\mu+2}(at))
    \right|,
\end{split}
\end{equation*}
where 
\begin{equation*}
    F_\mu(r,t)=\frac{e^{-\mu\phi_\mu(r)}e^{(\mu+2)\phi_{\mu+2}(t)}}{(\mu^2-r^2)^{\frac{1}{4}}((\mu+2)^2-t^2)^{\frac{1}{4}}}.
\end{equation*}
Using \eqref{phi} and \eqref{alpha}, we can write
\begin{equation*}
    F_\mu(r,t)=r^\mu t^{-\mu-2}f_\mu(r,t),
\end{equation*}
with
\begin{equation*}
\begin{split}
    f_\mu(r,t)=&
    \left(\frac{\mu+2+\sqrt{(\mu+2)^2-t^2}}{\mu+\sqrt{\mu^2-r^2}}\right)^{\mu}
    e^{\sqrt{\mu^2-r^2}-\sqrt{(\mu+2)^2-t^2}}
    \\
    &\quad\times
    \frac{\left(\mu+2+\sqrt{(\mu+2)^2-t^2}\right)^2}{(\mu^2-r^2)^{\frac{1}{4}}((\mu+2)^2-t^2)^{\frac{1}{4}}}.
\end{split}
\end{equation*}
Taking into account that 
\begin{equation*}
    a^2F_\mu(ar,at)=r^\mu t^{-\mu-2}f_\mu(ar,at),
\end{equation*}
and using the triangle inequality, we have that
\begin{equation*}
    D_{\mu,a}^2(r,t)\lesssim \sum_{m=1}^{4}D_{\mu,a}^{2,m},
\end{equation*}
where
\begin{align}
\label{D21}
     D_{\mu,a}^{2,1}
     &=
     r^\mu t^{-\mu-2}\left|f_\mu(r,t)-f_\mu(ar,at)\right|,
     \\
\label{D22}
     D_{\mu,a}^{2,2}
     &=
     r^\mu t^{-\mu-2}\left|f_\mu(r,t)e_{1,\mu}(r)-f_\mu(ar,at)e_{1,\mu}(ar)\right|,
     \\
\nonumber
     D_{\mu,a}^{2,3}
     &=
     r^\mu t^{-\mu-2}\left|f_\mu(r,t)e_{2,\mu}(t)-f_\mu(ar,at)e_{2,\mu}(at)\right|,
     \\
\nonumber
     D_{\mu,a}^{2,4}
     &=
     r^\mu t^{-\mu-2}\left|f_\mu(r,t)e_{1,\mu}(r)e_{2,\mu}(t)-f_\mu(ar,at)e_{1,\mu}(ar)e_{2,\mu}(at)\right|.
\end{align}
Therefore, it is enough to prove the estimate for 
\begin{equation}
\label{I2m}
    \mathrm{I}_2^m=\int_{1-\frac{1}{2a}}^{\frac{\mu}{2a}}\left(\int_r^{r+\frac{1}{a}}
        \left|D_{\mu,a}^{2,m}(r,t)\right|^2
        V(t)t\,dt\right)V(r)r\,dr,
        \qquad m=1,2,3,4.
\end{equation}
To obtain the result for $\mathrm{I}_2^1$, it is necessary to take advantage of cancellation.
Using the mean value theorem we have that there exists $s\in(0,1)$ such that
\begin{equation*}
    \left|f_\mu(r,t)-f_\mu(ar,at)\right|=(a-1)
    \left|r\partial_rf_\mu(r_{a,s},t_{a,s})
    +t\partial_tf_\mu(r_{a,s},t_{a,s})\right|,
\end{equation*}
where
$r_{a,s}=((1-s)a+s)r$ and $t_{a,s}=((1-s)a+s)t$.
One can check that 
\begin{equation*}
\begin{split}
    \partial_rf_\mu(r,t)=&
    \frac{r}{2}
    \left(\frac{1}{\mu^2-r^2}-\frac{2}{\mu+\sqrt{\mu^2-r^2}}\right)
    \,f_\mu(r,t),
    \\
    \partial_tf_\mu(r,t)=&
    \frac{t}{2}
    \left(\frac{1}{(\mu+2)^2-t^2}+\frac{2}{\mu+2+\sqrt{(\mu+2)^2-t^2}}\right)
    \,f_\mu(r,t).
\end{split}
\end{equation*}
And also that
\begin{equation}
\label{cota_fmu}
    |f_\mu(r,t)|\lesssim \mu, \qquad 0<r<\frac{\mu}{2},\, r<t<r+1.
\end{equation}
From here, if $0<r<\frac{\mu}{2a}$ and $r<t<r+\frac{1}{a}$, we have that
$0<r_{a,s}<\frac{\mu}{2}$ and $r_{a,s}<t_{a,s}<r_{a,s}+1$, and therefore $|f_\mu(r_{a,s},t_{a,s})|\lesssim \mu$ and
\begin{equation*}
    \mu^2-r_{a,s}^2\gtrsim \mu^2,\qquad
    (\mu+2)^2-t_{a,s}^2\gtrsim \mu^2.
\end{equation*}
Taking this into account and using the expressions for $r_{a,s}$ and $t_{a,s}$, we obtain that
\begin{equation*}
    \left|f_\mu(r,t)-f_\mu(ar,at)\right|\lesssim
    \mu
    \left(\frac{r^2+t^2}{\mu^2}
    +\left|\frac{t^2}{\mu+2+\sqrt{(\mu+2)^2-t_{a,s}^2}}
    -\frac{r^2}{\mu+\sqrt{\mu^2-r_{a,s}^2}}\right|\,\right).
\end{equation*}
And since $0<t-r<1/a<1$ and $r/\mu<1/2a<1$, we get
\begin{equation*}
     \left|f_\mu(r,t)-f_\mu(ar,at)\right|\lesssim
     \mu
    \left(\frac{r^2}{\mu^2}+\frac{r}{\mu}\right)
    \lesssim
    r.
\end{equation*}
Using this estimate in \eqref{D21} we obtain 
\begin{equation*}
    \left|D_{\mu,a}^{2,1}(r,t)\right|\lesssim r^{\mu+1}t^{-\mu-2},
\end{equation*}
and from \eqref{I2m} we see that
\begin{align}
    \mathrm{I}_2^1
    &\lesssim
    \int_{1-\frac{1}{2a}}^{\frac{\mu}{2a}}\left(\int_r^{r+\frac{1}{a}} t^{-2\mu-4}V(t)t\,dt\right)V(r)r^{2\mu+3}\,dr
    \nonumber
    \\
    &\le 
    \int_{1-\frac{1}{2a}}^{\frac{\mu}{2a}}\left(\int_r^{r+\frac{1}{a}}
    \frac{tV(t)}{(t^2-r^2)^\frac{1}{2}}(t^2-r^2)^\frac{1}{2}t^{-2\mu-4}\,dt\right)V(r)r^{2\mu+3}\,dr.
    \label{I21truco}
\end{align}
For a fixed $r>0$, we introduce the function $A_{\mu,r}(t)=(t^2-r^2)^\frac{1}{2}t^{-2\mu-4}$ for $t>r$. One can check that this function reaches its absolute maximum at 
$t_r=r\,\sqrt{\frac{2\mu+4}{2\mu+3}}$, and 
for $0<r<\frac{\mu}{2a}$ we have that
$t_r\in\left(r,r+\frac{1}{a}\right)$. 
Thus
\begin{equation*}
    A_{\mu,r}(t)\le A_{\mu,r}\left(t_r\right)
    \le \mu^{-\frac{1}{2}}r^{-2\mu-3},
    \quad
    r<t<r+\frac{1}{a},0<r<\frac{\mu}{2a}.
\end{equation*}
Using this and \eqref{normaTakeuchi} in \eqref{I21truco} we have that
\begin{align}
    \mathrm{I}_2^1
    &\lesssim  
    \mu^{-\frac{1}{2}} |||V|||
    \int_{1-\frac{1}{2a}}^{\frac{\mu}{2a}}
    V(r)\,dr
\label{est_I21}
    \\
    &\le \mu^{-\frac{1}{2}} |||V|||^2.
\nonumber
\end{align}
To get the estimate for $\mathrm{I}_2^m$ with $m=2,3,4$, it is not necessary to take advantage of cancellation. For brevity, we will only see how to obtain the estimate for $\mathrm{I}_2^2$. By arguing in the same way, we would obtain the estimates for $\mathrm{I}_2^3$ and $\mathrm{I}_2^4$.

Using the triangle inequality in \eqref{D22} 
and making a change of variables in \eqref{I2m}, since $a\ge 1$, we can write that
\begin{equation*}
    \begin{split}
    \mathrm{I}_2^2
    \le&\int_{1-\frac{1}{2a}}^{\frac{\mu}{2a}}
    \left(\int_r^{r+\frac{1}{a}}
    |f_\mu(r,t)|^2
    |e_{1,\mu}(r)|^2
    t^{-2\mu-4}V(t)t\,dt
    \right)V(r)r^{2\mu+1}\,dr
    \\
    &+
    \int_{a-\frac{1}{2}}^{\frac{\mu}{2}}
    \left(\int_r^{r+1}
    |f_\mu(r,t)|^2
    |e_{1,\mu}(r)|^2
    t^{-2\mu-4}V_a(t)t\,dt
    \right)V_a(r)r^{2\mu+1}\,dr,
    \end{split}
\end{equation*}
where $V_a$ is defined by \eqref{V_a}. Since $a\ge 1$ we can use here \eqref{cota_fmu} and \eqref{errores_cota1} to write
\begin{equation*}
    \begin{split}
    \mathrm{I}_2^2
    \lesssim&\int_{1-\frac{1}{2a}}^{\frac{\mu}{2a}}
    \left(\int_r^{r+\frac{1}{a}}
    t^{-2\mu-4}V(t)t\,dt
    \right)V(r)r^{2\mu+1}\,dr
    \\
    &+
    \int_{a-\frac{1}{2}}^{\frac{\mu}{2}}
    \left(\int_r^{r+1}
    t^{-2\mu-4}V_a(t)t\,dt
    \right)V_a(r)r^{2\mu+1}\,dr,
    \end{split}
\end{equation*}
By arguing as we did to obtain \eqref{est_I21}, since $a\ge 1$, we see that
\begin{equation*}
\begin{split}
     \mathrm{I}_2^2
    &\lesssim  
    \mu^{-\frac{1}{2}} |||V|||
    \left(
    \int_{1-\frac{1}{2a}}^{\frac{\mu}{2a}}
    V(r)r^{-2}\,dr
    +
    \int_{a-\frac{1}{2}}^{\frac{\mu}{2}}
    V_a(r)r^{-2}\,dr
    \right)
    \\
    &\le
    \mu^{-\frac{1}{2}} |||V|||
    \int_{\frac{1}{2}}^{\frac{\mu}{2}}
    (V(r)+V_a(r))r^{-2}\,dr
    \\
    &\lesssim
    \mu^{-\frac{1}{2}} |||V|||^2,
\end{split}
\end{equation*}
In the last inequality we used the identity \eqref{normaV_a}.
This completes the proof of the lemma.
\hfill $\square$

The following lemma will be needed in the proof of Lemma \ref{escalones}, and its proof can be found in \cite[p. 446]{Watson}.
\begin{lemma}
    \label{decrecimiento}
    If $\mu>1/2$ then, $t|H^{(1)}_\mu(t)|^2$ is a decreasing function of $t$ in the interval $(0,\infty).$
\end{lemma}
\textit{Proof of Lemma \ref{escalones}.}
Fix $m\in\{0,\ldots,N\}$ and consider 
\begin{equation*}
   0<r<1+\frac{m-1}{2a} 
   \quad\text{ and }\quad
   t>1+\frac{m}{2a}.
\end{equation*}
Since $0<r<\mu/2$, from \eqref{cota_J1} we have 
\begin{equation*}
    \left|J_\mu(r)\right|^2\lesssim 
        \mu^{-1}e^{-2\mu\phi_\mu(r)}.
\end{equation*}
From \eqref{phi} one can see that
\begin{equation}
    \label{derivada_phi}
    \phi'(r)=-\frac{\sqrt{\mu^2-r^2}}{\mu r},
\end{equation}
so $\phi$ is a decreasing function in $(0,\infty)$, 
and therefore, we have 
\begin{equation}
\label{escalon1}
    r\left|J_\mu(r)\right|^2\lesssim 
        \left(1+\frac{m-1}{2a}\right)\mu^{-1}
        e^{-2\mu\phi_\mu\left(1+\frac{m-1}{2a}\right)},
        \qquad
        0<r<1+\frac{m-1}{2a}.
\end{equation}
On the other hand, from Lemma \ref{decrecimiento} we know that $t|H^{(1)}_\mu(t)|^2$ is a decreasing function in the interval $(0,\infty),$ and therefore
\begin{equation*}
    t\left|H^{(1)}_{\mu+2}(t)\right|^2\le \left(1+\frac{m}{2a}\right) \left|H^{(1)}_{\mu+2}\left(1+\frac{m}{2a}\right)\right|^2,
    \qquad
    t>1+\frac{m}{2a}.
\end{equation*}
From here, since $1+\frac{m}{2a}<\frac{\mu+2}{2}$, using \eqref{cota_H1}, we get
\begin{equation}
\label{escalon2}
    t\left|H^{(1)}_{\mu+2}(t)\right|^2\lesssim 
    \left(1+\frac{m}{2a}\right)\mu^{-1}e^{2(\mu+2)\phi_{\mu+2}\left(1+\frac{m}{2a}\right)},
    \qquad
    t>1+\frac{m}{2a}.
\end{equation}
From \eqref{escalon1} and \eqref{escalon2} we obtain that for $m=0,\ldots,N-1,$
\begin{equation*}
\begin{split}
    \mathrm{L}_m&\lesssim
    h_\mu(m)
    \int_0^\infty V(r)\,dr
    \int_{1+\frac{m}{2a}}^{1+\frac{m+1}{2a}}V(t)\,dt,\\
    \mathrm{L}_{N}&\lesssim
    h_\mu(N)
    \int_0^\infty V(r)\,dr
    \int_{1+\frac{N}{2a}}^{\infty}V(t)\,dt,
\end{split}
\end{equation*}
where
\begin{equation*}
    h_\mu(m)=
    \left(1+\frac{m-1}{2a}\right)
    \left(1+\frac{m}{2a}\right)
    \mu^{-2}
    e^{-2\mu\phi_\mu\left(1+\frac{m-1}{2a}\right)}
    e^{2(\mu+2)\phi_{\mu+2}\left(1+\frac{m}{2a}\right)}.
\end{equation*}
Therefore,
\begin{equation*}
\begin{split}
    \sum_{m=0}^{N}\mathrm{L}_m&\lesssim
    \max_{0 \leq m \leq N} h_\mu(m) 
    \int_0^\infty V(r)\,dr
    \left(
    \sum_{m=0}^{N-1}
    \int_{1+\frac{m}{2a}}^{1+\frac{m+1}{2a}}V(t)\,dt
    +
    \int_{1+\frac{N}{2a}}^{\infty}V(t)\,dt
    \right)
    \\
    &\le
     \max_{0 \leq m \leq N} h_\mu(m)\, 
    |||V|||^2.   
\end{split}
\end{equation*}
Thus, it is enough to prove that for $\mu$ large enough, 
\begin{equation*}
    \max_{0 \leq m \leq N} h_\mu(m)\lesssim 1.
\end{equation*}
Using the identities \eqref{alpha} and \eqref{phi} we write
\begin{multline*}
    h_\mu(m)=
    \left(1+\frac{m-1}{2a}\right)
    \left(1+\frac{m}{2a}\right)
    \mu^{-2}
    e^{-2\left(\sqrt{(\mu+2)^2-\left(1+\frac{m}{2a}\right)^2}-\sqrt{\mu^2-\left(1+\frac{m-1}{2a}\right)^2}\right)}\\
    \times
    \left(\frac{\mu+\sqrt{\mu^2-\left(1+\frac{m-1}{2a}\right)^2}}{1+\frac{m-1}{2a}}\right)^{-2\mu}
    \left(\frac{\mu+2+\sqrt{(\mu+2)^2-\left(1+\frac{m}{2a}\right)^2}}{1+\frac{m}{2a}}\right)^{2\mu+4}.
\end{multline*}
One can check that for fixed $m\in\{0,1,\ldots,[\mu-2a]\}$,
\begin{equation*}
    h_\mu(m)\le
    \left(\frac{1+\frac{m-1}{2a}}{1+\frac{m}{2a}}\right)^{2\mu+1}
    \frac{(2(\mu+2))^4}{\mu^2\left(1+\frac{m}{2a}\right)^2}
    \left(\frac{\mu+2+\sqrt{(\mu+2)^2-\left(1+\frac{m}{2a}\right)^2}}{\mu+\sqrt{\mu^2-\left(1+\frac{m-1}{2a}\right)^2}}\right)^{2\mu}.
\end{equation*}

Now we consider $x\in\left(0,\frac{\mu}{2a}-1\right)$ and introduce the functions
\begin{equation*}
    \begin{split}
    f_\mu(x)&= \left(1-\frac{1}{2a(1+x)}\right)^{2\mu+1}
    \frac{(2(\mu+2))^4}{\mu^2\left(1+x\right)^2} \ \ {\rm and}
    \\
    g_\mu(x)&=\left(\frac{\mu+2+\sqrt{(\mu+2)^2-\left(1+x\right)^2}}{\mu+\sqrt{\mu^2-\left(1+x\right)^2}}\right)^{2\mu}.
    \end{split}
\end{equation*}
One can check that for a fixed $m\in\{0,1,\ldots,[\mu-2a]\}$,
\begin{equation*}
    h_\mu(m)\le f_\mu\left(\frac{m}{2a}\right)g_\mu\left(\frac{m}{2a}\right).
\end{equation*}
Since $f_\mu(x)$ and $g_\mu(x)$ are increasing functions in the interval $\left(0,\frac{\mu}{2a}-1\right)$, we have that for all $m\in\{0,1,\ldots,[\mu-2a]\}$,
\begin{equation*}
\begin{split}
     h_\mu(m)
     &\le  
     f_\mu\left(\frac{\mu}{2a}-1\right)g_\mu\left(\frac{\mu}{2a}-1\right)
     \\
     &\lesssim
     \left(\frac{\mu+2+\sqrt{(\mu+2)^2-\frac{\mu^2}{4a}}}{\mu+\sqrt{\mu^2-\frac{\mu^2}{4a}}}\right)^{2\mu}.
\end{split}
\end{equation*}
From here,
\begin{equation}
\label{h}
    \lim_{\mu\rightarrow+\infty}h_\mu(m)\lesssim 1
\end{equation}
uniformly in $m$ with $m\in\{0,1,\ldots,[\mu-2a]\}$.
\hfill
$\square$

To prove Lemma \ref{Juan}, we will follow the strategy used in \cite[pp. 370]{BRV1997} to bound $\mathscr{I}_3$.

\bigskip
\textit{Proof of Lemma \ref{Juan}.}
For fixed $a\ge1$, we denote by $\mathrm{I}$ the iterated integral in \eqref{est_Juan}.
Since $H_{\mu+2}=J_{\mu+2}+iY_{\mu+2}$, using the triangle inequality, we split $\mathrm{I}$ into two
\begin{equation*}
    \mathrm{I}\lesssim \mathrm{I}_1+\mathrm{I}_2,
\end{equation*}
where
\begin{equation*}
    \begin{split}
        \mathrm{I}_1
        =&
        \int_{\frac{\mu}{2a}-\frac{1}{a}}^{\mu+2-(\mu+2)^{\frac{1}{3}}}\int_r^{\mu+2-(\mu+2)^{\frac{1}{3}}}
        |J_{\mu+2}(t)|^2
        V(t)t\,dt\,|J_\mu(r)|^2V(r)r\,dr,
        \\
        \mathrm{I}_2
        =&
        \int_{\frac{\mu}{2a}-\frac{1}{a}}^{\mu+2-(\mu+2)^{\frac{1}{3}}}\int_r^{\mu+2-(\mu+2)^{\frac{1}{3}}}
        |Y_{\mu+2}(t)|^2
        V(t)t\,dt\,|J_\mu(r)|^2V(r)r\,dr.
    \end{split}
\end{equation*}
Using Lemma \ref{lema_BRV}, we estimate $\mathrm{I}_1$:
\begin{equation*}
\begin{split}
    \mathrm{I}_1
    &\le 
    \int_0^\infty |J_{\mu}(r)|^2 V(r)r\,dr
    \int_0^\infty |J_{\mu+2}(t)|^2 V(t)t\,dt
    \\
    &
    \lesssim |||V|||^2.
\end{split}
\end{equation*}
To obtain the estimate for $\mathrm{I}_2$, we use the Debye expansions given in \eqref{Deby_J} and \eqref{Deby_Y} and the bounds for the error terms given in \eqref{errores_cota2} to write
\begin{equation}
\label{I2_juan}
   \mathrm{I}_2\lesssim 
   \int_{\frac{\mu}{2a}-\frac{1}{a}}^{\mu+2-(\mu+2)^{\frac{1}{3}}}
   \int_r^{\mu+2-(\mu+2)^{\frac{1}{3}}}
   \frac{e^{2(\mu+2)\phi_{\mu+2}(t)}}{((\mu+2)^2-t^2)^{\frac{1}{2}}}
   V(t)t\,dt\,
   \frac{e^{-2\mu\phi_{\mu}(r)}}{(\mu^2-r^2)^{\frac{1}{2}}}V(r)r\,dr.
\end{equation}
From \eqref{phi} and \eqref{alpha}, one can verify that for $r$ in the integration region, if $\mu>\mu_0$ with $\mu_0=\mu_0(a)$ large enough, then
\begin{equation}
\label{control_exponencial}
    e^{2(\mu+2)\phi_{\mu+2}(r)-2\mu\phi_{\mu}(r)}\sim 1.
\end{equation}
Using this in the previous inequality, 
we obtain 
\begin{equation*}
\begin{split}
    \mathrm{I}_2&\lesssim 
   \int_{\frac{\mu}{2a}-\frac{1}{a}}^{\mu+2-(\mu+2)^{\frac{1}{3}}}
   \int_r^{\mu+2-(\mu+2)^{\frac{1}{3}}}
   \frac{e^{2(\mu+2)(\phi_{\mu+2}(t)-\phi_{\mu+2}(r))}}{((\mu+2)^2-t^2)^{\frac{1}{2}}(\mu^2-r^2)^{\frac{1}{2}}}
   V(t)t\,dt\,
   V(r)r\,dr
   \\
   &=
   \int_{\frac{\mu}{2a}-\frac{1}{a}}^{\mu+2-(\mu+2)^{\frac{1}{3}}}
   \int_r^{\mu+2-(\mu+2)^{\frac{1}{3}}}
   \frac{e^{2(\mu+2)(t-r)\phi'_{\mu+2}(b_{r,t})}}{((\mu+2)^2-t^2)^{\frac{1}{2}}(\mu^2-r^2)^{\frac{1}{2}}}
   V(t)t\,dt\,
   V(r)r\,dr,
\end{split}
\end{equation*}
where $r<b_{r,t}<t.$
Now we decompose the integration interval in $r$ dyadically as follows:
\begin{equation*}
    \left(\frac{\mu}{2a}-\frac{1}{a},\mu+2-(\mu+2)^{\frac{1}{3}}\right)
    \subseteq
    \bigcup_{j=0}^N A_j,    
\end{equation*}
where 
\begin{equation*}
    A_j=\left[\mu+2-2^{j+1}(\mu+2)^{\frac{1}{3}},\mu+2-2^{j}(\mu+2)^{\frac{1}{3}}\right),
    \qquad j=0,\ldots,N,
\end{equation*} 
and $N$ is the smallest natural number satisfying 
\begin{equation*}
    2^{N+1}\ge\left(\mu+2-\frac{\mu}{2a}+\frac{1}{a}\right)(\mu+2)^{-\frac{1}{3}}.
\end{equation*}
From here, we have 
\begin{equation}
\label{Juan_sumaIj}
    \mathrm{I}_2\le \sum_{j=0}^{N}\mathrm{L}_j,
\end{equation}
where
\begin{equation*}
    \mathrm{L}_j=
    \int_{A_j}
   \int_{B_r}
   \frac{e^{2(\mu+2)(t-r)\phi'_{\mu+2}(b_{r,t})}}{((\mu+2)^2-t^2)^{\frac{1}{2}}(\mu^2-r^2)^{\frac{1}{2}}}
   V(t)t\,dt\,V(r)r\,dr,
   \qquad j=0,\ldots N,
\end{equation*}
and $B_r=\left(r,\mu+2-(\mu+2)^{\frac{1}{3}}\right)$.
Now introducing the function
\begin{equation*}
    g_j(r,t)=f(r,t)
    \left(r^2-\left(\mu+2-2^{j+1}(\mu+2)^{\frac{1}{3}}\right)^2\right)^{\frac{1}{2}}
    (t^2-r^2)^{\frac{1}{2}},
\end{equation*}
where
\begin{equation*}
    f(r,t)=\frac{e^{2(\mu+2)(t-r)\phi'_{\mu+2}(b_{r,t})}}{((\mu+2)^2-t^2)^{\frac{1}{2}}(\mu^2-r^2)^{\frac{1}{2}}},
\end{equation*}
we can write
\begin{align}
    \mathrm{L}_j&\le
    \underset{\begin{subarray}{l}t\in B_r\\r\in A_j\end{subarray}}{\text{sup}}g_j(r,t)
    \int_{A_j}
   \int_{B_r}
   \frac{tV(t)}{(t^2-r^2)^{\frac{1}{2}}}\,dt\,
   \frac{rV(r)}{\left(r^2-\left(\mu+2-2^{j+1}(\mu+2)^{\frac{1}{3}}\right)^2\right)^{\frac{1}{2}}}\,dr
   \nonumber
   \\
   &\le
   \underset{\begin{subarray}{l}t\in B_r\\r\in A_j\end{subarray}}{\text{sup}}g_j(r,t)\,
   |||V|||^2.
   \label{Juan_Ij}
\end{align}
One can check that
\begin{equation}
\label{supremo_gj}
    \underset{\begin{subarray}{l}t\in B_r\\r\in A_j\end{subarray}}{\text{sup}}g_j(r,t)
    \lesssim
    \underset{\begin{subarray}{l}t\in B_r\\r\in A_j\end{subarray}}{\text{sup}}
    \frac{(t-r)^{\frac{1}{2}}
    e^{2(\mu+2)(t-r)\phi'_{\mu+2}(b_{r,t})}}
    {(\mu+2-t)^{\frac{1}{2}}}.
\end{equation}
From \eqref{derivada_phi} one can see that 
$\phi'_\mu$ is an increasing function, and since $r<b_{r,t}<t$, we have that
\begin{equation*}
    \phi'_{\mu+2}(b_{r,t})\le\phi'_{\mu+2}(t)=-\frac{((\mu+2)^2-t^2)^{\frac{1}{2}}}{(\mu+2) t}.
\end{equation*}
Using this in \eqref{supremo_gj}, we get
\begin{equation}
\label{supremo}
    \underset{\begin{subarray}{l}t\in B_r\\r\in A_j\end{subarray}}{\text{sup}}g_j(r,t)
    \lesssim
    \underset{\begin{subarray}{l}t\in B_r\\r\in A_j\end{subarray}}{\text{sup}}
    \frac{(t-r)^{\frac{1}{2}}
    e^{-\mu^{-\frac{1}{2}}(t-r)(\mu+2-t)^{\frac{1}{2}}}}
    {(\mu+2-t)^{\frac{1}{2}}}.
\end{equation}
To calculate the supremum, 
for a fixed $r\in A_j$, we decompose $B_r$ dyadically as follows:
\begin{equation}
\label{Br}
    B_r
    =
    B_r^1\cup B_r^2,
\end{equation}
with
\begin{equation}
\begin{array}{lll}
    B_r^1=A_0\cup\ldots\cup A_{j-2},
    \qquad
    &
    B_r^2\subseteq A_j\cup A_{j-1},\qquad 
    &
    j\ge 2,
    \\
    B_r^1=\emptyset,
    \qquad
    &
    B_r^2\subseteq A_1\cup A_0,\qquad 
    &j=1,
    \\
    B_r^1=\emptyset,
    \qquad
    &
    B_r^2\subseteq A_0,\qquad 
    &j=0.
\end{array}
\label{Br12}
\end{equation}
Observe that if 
$r\in A_j$ and 
$t\in B_r^1$, then 
$t-r\sim 2^j\mu^{\frac{1}{3}}$
and
$\mu+2-t\gtrsim \mu^{\frac{1}{3}}$. Thus
\begin{equation}
\label{supremo1}
    \underset{\begin{subarray}{l}t\in B_r^1\\r\in A_j\end{subarray}}{\text{sup}}
    \frac{(t-r)^{\frac{1}{2}}
    e^{-\mu^{-\frac{1}{2}}(t-r)(\mu+2-t)^{\frac{1}{2}}}}
    {(\mu+2-t)^{\frac{1}{2}}}
    \lesssim
    2^{\frac{j}{2}}
    e^{-2^{\frac{j}{2}}}.
\end{equation}
On the other hand, 
if $r\in A_j$ and $t\in B_r^2$, then 
$0\le t-r\lesssim 2^j\mu^{\frac{1}{3}}$
and
$\mu+2-t\gtrsim 2^j\mu^{\frac{1}{3}}$. Hence
\begin{align}
\nonumber
    \underset{\begin{subarray}{l}t\in B_r^2\\r\in A_j\end{subarray}}{\text{sup}}
    \frac{(t-r)^{\frac{1}{2}}
    e^{-\mu^{-\frac{1}{2}}(t-r)(\mu+2-t)^{\frac{1}{2}}}}
    {(\mu+2-t)^{\frac{1}{2}}}
    &\lesssim
    \underset{0<x<2^{j}\mu^{\frac{1}{3}}}{\text{sup}}
    \frac{x^{\frac{1}{2}}
    e^{-\mu^{-\frac{1}{2}}x2^{\frac{j}{2}}\mu^{\frac{1}{6}}}}
    {2^{\frac{j}{2}}\mu^{\frac{1}{6}}}
    \\
\nonumber
    &\le
    \underset{0<y<2^{\frac{3j}{2}}}{\text{sup}}
    2^{-\frac{3j}{4}}y^{\frac{1}{2}}
    e^{-y}
    \\
    &\le
    2^{-\frac{3j}{4}}.
\label{supremo2}
\end{align}
Using \eqref{Br}-\eqref{supremo2} in \eqref{supremo}
we get
\begin{equation*}
    \begin{split}
    \underset{\begin{subarray}{l}t\in B_r\\r\in A_j\end{subarray}}{\text{sup}}g_j(r,t)&\lesssim
    2^{-\frac{3j}{4}},
    \end{split}
\end{equation*}
and inserting this inequality in \eqref{Juan_Ij} we obtain
\begin{equation*}
    \mathrm{L}_j\lesssim 2^{-\frac{3j}{4}}|||V|||^2.
\end{equation*}
Inserting this inequality into \eqref{Juan_sumaIj} concludes the proof of the lemma.
\hfill$\square$

\bigskip
\textit{Proof of Lemma \ref{bueno}.}
The proof of this lemma is almost the same as that of Lemma \ref{Juan}, but in this case, instead of \eqref{I2_juan} we have
\begin{equation*}
   \mathrm{I}_2\lesssim 
   \int_0^{\mu-\mu^{\frac{1}{3}}}
   \int_r^{\mu-\mu^{\frac{1}{3}}}
   \frac{e^{2\mu\phi_{\mu}(t)}}{(\mu^2-t^2)^{\frac{1}{2}}}
   V(t)t\,dt\,
   \frac{e^{-2(\mu+2)\phi_{\mu+2}(r)}}{((\mu+2)^2-r^2)^{\frac{1}{2}}}
   V(r)r\,dr.
\end{equation*}
So instead of \eqref{control_exponencial}, we use that for $r$ in the integration region, if $\mu$ is large enough, then
\begin{equation}
\label{control_exponencial_bueno}
    e^{2\mu\phi_{\mu}(r)-2(\mu+2)\phi_{\mu+2}(r)}\lesssim 1,
\end{equation}
and hence
\begin{equation*}
   \mathrm{I}_2\lesssim 
   \int_0^{\mu-\mu^{\frac{1}{3}}}
   \int_r^{\mu-\mu^{\frac{1}{3}}}
   \frac{e^{2\mu(\phi_{\mu}(t)-\phi_{\mu}(r))}}{(\mu^2-t^2)^{\frac{1}{2}}((\mu+2)^2-r^2)^{\frac{1}{2}}}
   V(t)t\,dt\,
   V(r)r\,dr.
\end{equation*}
To complete the proof of the lemma, we argue as in Lemma \ref{Juan}. 
\hfill$\square$
\begin{remark}
    \label{explica_bueno}
    We would like to point out that in the estimate \eqref{est_buena}, the lower limit of the integral in $r$ is $0$, while in the estimate \eqref{est_Juan}, we cannot reach this value. This is because \eqref{control_exponencial_bueno} is true for small values of $r$, while \eqref{control_exponencial} is not.
\end{remark}

\end{document}